\documentclass[11pt,a4paper]{scrarticle}
\usepackage[utf8]{inputenc}
\usepackage[T1]{fontenc}
\usepackage[USenglish]{babel}
\usepackage[singlespacing]{setspace}
\usepackage{paralist}

\usepackage{cite}
\usepackage{graphicx}
\usepackage{verbatim}
\usepackage{amsmath}
\usepackage{amsfonts}
\usepackage{amssymb}
\usepackage{amsthm}
\usepackage{rotating}
\usepackage{upgreek}
\usepackage{braket}
\usepackage{hyperref}
\usepackage{mathtools}
\usepackage{marginnote}
\usepackage{MnSymbol}
\usepackage{chemfig}
\usepackage{chemformula}
\usepackage[version=3]{mhchem}
\usepackage{tikz}
\usepackage{units}
\usepackage{layout}   
\usepackage{float}
\usepackage{subfigure}
\usepackage{caption}
\usepackage{hyperref}
\usepackage{multibib}
\usepackage{makeidx}
\usepackage{bbm}
\makeindex

\author{Christian Streib}
\title{Minimal Presentation of $PSL(2,\Z)$ Using Continuant Matrices with Integer Coefficients}
\date{\today}

\newcommand{\C}{\mathbb{C}}

\newcommand{\N}{\mathbb{N}}
\newcommand{\Z}{\mathbb{Z}}

\swapnumbers
\theoremstyle{definition}
\newtheorem{df}{Definition}[section]
\newtheorem{sz}[df]{Theorem}
\newtheorem{genericthm}[df]{\thistheoremname}
\newtheorem{bsp}[df]{Example}
\newtheorem{cor}[df]{Corollary}
\newtheorem{rem}[df]{Remark}
\newtheorem{lm}[df]{Lemma}
\newenvironment{bw}[1][\textbf{Proof}]{\begin{proof}[#1]}{\end{proof}}

\newcommand{\thistheoremname}{}
\newenvironment{bla}[1]
{\renewcommand{\thistheoremname}{#1}\begin{genericthm}}{\end{genericthm}}
\numberwithin{equation}{section}

\begin{document}

\begin{center}
{\huge\bfseries Minimal Presentation of $PSL(2,\Z)$ \par Using Continuant Matrices with \par Integer Coefficients \par}
\vspace{0.5cm}
by\par
{\Large Christian Helmut Anton \textsc{Streib}\par}
\vspace{0.5cm}
{\large \today \par}
Stuttgart
\end{center}

\section*{Abstract}

In this article, the goal is to find the shortest presentation of a matrix $A \in PSL(2,\Z)$ in terms of the so-called continuant matrices which are most known for their role in continued fraction theory. In chapter 7 of \cite{Farey}, Morier-Génoud and Ovsienko investigate this problem with the restriction that all coefficients of the continuant matrices are positive. Now, the goal is to determine the shortest presentation allowing all integer coefficients.
To determine this minimal presentation, a few characteristic transformations will be introduced. It will also be investigated under which conditions such a minimal presentation becomes unique. The results are also generalized on conjugacy classes in $PSL(2,\Z)$.\\
If you wish to contact the author or you have some questions related to this work, feel free to write an email to \href{mailto:christian.streib@online.de}{christian.streib@online.de}.
\tableofcontents
\setlength{\parindent}{0mm}

\section{Introduction}
\label{ch:intro}

The idea of this article is to determine the minimal presentation of a matrix $A \in PSL(2,\Z)$ as a product of single continuant matrices, i.e. matrices of the form $\begin{pmatrix}c&-1\\1&0\end{pmatrix}$ where $c$ can take any integer values. Concretely, the goal is to write $A=M(c_1,...,c_n):=M(c_1) \cdot ... \cdot M(c_n)$. In general, for a matrix in $PSL(2,\Z)$, there are many different possibilities to express it in such a way. Now, for a given $A \in PSL(2,\Z)$, the shortest sequence $(c_1,...,c_n)$ such that this condition is fulfilled is searched, i.e. the number $n$ shall be minimized. The crucial parts of this work are the transformations in lemma \ref{lm:surg}, their inverse operations (lemma \ref{lm:inverse}), the transformations introduced under numeral \ref{bla:newsurgery} and in lemma \ref{lm:232} and theorem \ref{sz:calc2} which allows in combination with theorem \ref{sz:sequence} to deduce that the transformations mentioned above are sufficient to obtain the minimal presentation.\\
The whole topic is also closely related to quiddity cycles, i.e. sequences $(c_1,...,c_n)$ such that $M(c_1,...,c_n)=Id$ in $PSL(2,\Z)$. For more details concerning them, see also \cite{Ovsienko} and \cite{Streib}. Quiddity cycles are especially closely related to conjugation in $PSL(2,\Z)$. Using some elementary properties of quiddity cycles, the results about minimality and uniqueness of a presentation of a matrix in $PSL(2,\Z)$ will also be extended on conjugacy classes.\\

In section \ref{sec:huit}, some basic terms for this article such as continuant matrices and quiddity cycles are defined. Furthermore, a few identities and transformations introduced in chapters 4 and 6 of \cite{Cuntz} are recalled. The principles of continued fraction expressions and negative continued fraction expressions of rationals are explained and a modification of the euclidean algorithm to obtain the negative continued fraction expression of a rational is introduced.\\
In section \ref{sec:un}, the central statement is theorem \ref{sz:calc2} stating that for each matrix in $PSL(2,\Z)$, there is a unique presentation with continuant matrices involving the negative continued fraction expression of the ratio of the matrix's two left entries. The presentation obtained there is very similar to the one given in \cite{Farey}, in some cases even identical (see corollary \ref{cor:misone}). However, allowing negative coefficients as well, there may still be a shorter presentation.\\
To find this shortest presentation is the main goal in section \ref{sec:deux}. To do so, more identities for continuant matrices and corresponding transformations of sequences are introduced. Using the fact that the presentation from the previous section (see theorem \ref{sz:calc2}) is unique, it can be shown that together with the transformations of \cite{Cuntz}, these transformations suffice to obtain the minimal presentation of a given matrix in $PSL(2,\Z)$. Furthermore, it turns out that this minimal presentation becomes unique if the value -2 is forbidden as a coefficient except for the first or for the last entry.\\
In section \ref{sec:conj}, these results are generalized on conjugacy classes. The central finding of this section is that basically all rules that allow to determine whether a presentation is minimal for a single matrix also hold for conjugacy classes, but, for single matrices, all these rules have exceptions for the first and the last coefficient of the sequence, and these exceptions fall away in the case of conjugacy classes.
\section{Basic Terms and Definitions}
\label{s:continuant}
\label{sec:huit}

\subsection{Continuant Matrices and Quiddity Cycles}

In this section, the basic terms of continuant matrices and quiddity cycles will be defined and explained. Furthermore, a few identities for continuant matrices are introduced as well as the so-called surgery operations related to them which are needed in section \ref{sec:deux}. In the first subsection, the central terms of continuant matrices and quiddity cycles are introduced. The topic of the second subsection are so-called surgery operations. In the third subsection, the concept of continued fraction expressions and negative continued fraction expressions is recalled. Furthermore, a variant of the euclidean algorithm to find the negative continued fraction expression of a given ratio is introduced.\\
For background information about the quiddity cycles and surgery operations, see also \cite{Cuntz} and \cite{Farey}. For details about continued fraction theory in general, see chapter 10 in \cite{Rosen}. More concrete explanations about negative continued fraction expressions and their relation to regular continued fraction expressions and triangulations, can be found in chapter 2 of \cite{Farey}.

\begin{df} \begin{em}
For any given $c \in \C$, consider the matrix
\begin{align}
M(c):=\begin{pmatrix}c&-1\\1&0\end{pmatrix}.
\end{align}
The continuant matrix or matrix of continued fractions as it is called in \cite{Farey} $M(c_1,...,c_n)$ is then defined as the product of the matrices $M(c_1),...,M(c_n)$, that is
\begin{align}
M(c_1,...,c_n):=M(c_1) \cdot ... \cdot M(c_n).
\label{eqn:continuant}
\end{align} \end{em}
\end{df}

\begin{bsp}
The continuant matrix $M(3,5,-2)$ is defined as $M(3) \cdot M(5) \cdot M(-2)$, and explicitly, one gets
\begin{equation}
A=\begin{pmatrix}3&-1\\1&0\end{pmatrix}\begin{pmatrix}5&-1\\1&0\end{pmatrix}\begin{pmatrix}-2&-1\\1&0\end{pmatrix}=\begin{pmatrix}14&-3\\5&-1\end{pmatrix}\begin{pmatrix}-2&-1\\1&0\end{pmatrix}=\begin{pmatrix}-31&-14\\-11&-5\end{pmatrix}.
\end{equation}
\end{bsp}

In the following, continuant matrices are used to define quiddity cycles, objects that play a crucial role in so-called Frieze patterns invented by Coxeter and Conway \cite{Coxeter}. However, here the quiddity cycles are merely used as a tool to simplify the argumentation to construct a unique minimal presentation for matrices in $PSL(2,\Z)$ and their conjugacy classes. The surgery operations introduced in lemma \ref{lm:surg} which can be applied on quiddity cycles and their inverse transformations out of lemma \ref{lm:inverse} play a crucial role for nearly all of the following chapters.

\begin{df} \begin{em}
If there is an $\varepsilon \in \{\pm 1 \}$ such that the equation
\begin{align}
M(c_1,...,c_n)=\varepsilon Id
\label{eqn:epsilon}
\end{align}
is fulfilled, a sequence $(c_1,...,c_n) \in \Z$ is called a quiddity cycle. The term cycle is used because of the fact that every cyclic permutation of a quiddity cycle is also a quiddity cycle. $\varepsilon$ will be called the sign of the quiddity cycle. $n$ is called the length of the quiddity cycle.
\label{df:qc} \end{em}
\end{df}

\begin{bsp}
One can calculate that $M(1,3,1,2,2)=M(1)M(3)M(1)M(2)M(2)=-Id$. Hence, $(1,3,1,2,2)$ is a quiddity cycle of sign $-1$.
\label{bsp:qc}
\end{bsp}

\begin{lm}
Let $R \subset \C$ be a ring and $a,b,u \in R$. Then, the following three equations hold:
\begin{align}
M(a)M(b)&=M(a+1)M(1)M(b+1)
\label{eqn:calc1}\\
M(a)M(b)&=-M(a-1)M(-1)M(b-1)
\label{eqn:calc2}\\
M(a)M(0)M(b)&=M(a+u)M(0)M(b-u)=-M(a+b)
\label{eqn:calc3}
\end{align}
The proof of this lemma is straightforward calculation that is left to the reader.
\label{lm:calc}
\end{lm}

\begin{lm}
Let $(c_1,...,c_n)$ be a quiddity cycle, i.e. $M(c_1,...,c_n)=\pm Id$. Then, any cyclic permutation of $(c_1,...,c_n)$ is also a quiddity cycle.
\label{lm:fullrot}
\end{lm}
\begin{bw}
 Since the identity matrix commutes with every other matrix in $SL(2,\C)$, for $1 < k \leq n$, it is
\begin{equation}
\begin{aligned}
&M(c_k,...,c_n,c_1,...,c_{k-1})\\
&~~=M^{-1}(c_{k-1}) \cdot ... \cdot M^{-1}(c_1) \cdot  M(c_1,...,c_n) \cdot M(c_1) \cdot ... \cdot M(c_{k-1}) \\
&~~=M^{-1}(c_{k-1}) \cdot ... \cdot M^{-1}(c_1) \cdot  (\pm Id) \cdot M(c_1) \cdot ... \cdot M(c_{k-1}) \\
&~~=\pm Id
\end{aligned}
\end{equation}
Thus, for a sequence $(c_1,...,c_n)$, the property of solving the equation $M(c_1,...,c_n)=\pm Id$ is invariant under cyclic permutations of its entries.
\end{bw}

\subsection{Surgery Operations}

\begin{bla}{Lemma: Surgery Operations}
Any quiddity cycle of length at least three containing $c_i \in \{0,\pm1\}$ can be reduced by one of the following three operations, and the result is another shorter quiddity cycle. To keep the notation consistent with \cite{Farey}, those operations are called surgery operations\index{surgery operations} in what follows.\\
(a) Case $c_i=1$:
\begin{align}
(c_1,...,c_{i-1},1,c_{i+1},...,c_n) \rightarrow (c_1,...,c_{i-1}-1,c_{i+1}-1,...,c_n)
\label{eqn:surgery1}
\end{align}
(b) Case $c_i=-1$:
\begin{align}
(c_1,...,c_{i-1},-1,c_{i+1},...,c_n) \rightarrow (c_1,...,c_{i-1}+1,c_{i+1}+1,...,c_n)
\label{eqn:surgery2}
\end{align}
(c) Case $c_i=0$:
\begin{align}
(c_1,...,c_{i-1},0,c_{i+1},...,c_n) \rightarrow (c_1,...,c_{i-1}+c_{i+1},...,c_n)
\label{eqn:surgery3}
\end{align}
Obviously, in the first two cases, the newly obtained shorter quiddity cycle is of length $n-1$, and in case (c), it is of length $n-2$. Furthermore, in case (a), the sign of the two quiddity cycles is the same, whilst in cases (b) and (c), the reduced one has the opposite sign compared to the initial one. Also note that these reductions still work if $i=1$ or $i=n$. If $i=1$, the entry $c_n$ then takes the role of $c_{i-1}$ with respect to the modifications, and if $i=n$, the entry $c_1$ takes the role of $c_{i+1}$ analogously. However, in these cases, the represented matrix $M(c_1,...,c_n)$ might change in general; for more details, see section \ref{sec:conj}.
\label{lm:surg}
\end{bla}
\begin{bw}
If $(c_1,...,c_n)$ is a quiddity cycle, one of the equations
\begin{align}
M(c_1,...,c_n)=\pm Id
\label{eqn:quiddity}
\end{align}
holds. Applying equation (\ref{eqn:continuant}) and lemma \ref{lm:calc}, one can see that the reduced expressions still solve one of the equations above. Hence, they are also quiddity cycles. The different sign of the reduced quiddity cycles in the cases (b) and (c) comes from the minus sign in the equations (\ref{eqn:calc2}) and (\ref{eqn:calc3}). In the cases $i=1$ and $i=n$, rotate the cycle in such a way that the entry $c_i$ is now an inner entry by applying the previous lemma for a suitable value of $k$. Then perform the reduction and revert the rotation to obtain the same result. The statement about the length of the reduced cycles is obvious.
\label{bla:surgery}
\end{bw}

\begin{bla}{Inverse Surgery Operations}
To make notations in the following chapters more understandable, here, the operations inverse to those of lemma \ref{lm:surg} are presented\index{inverse surgery operations}:\\
(a) Inserting a one and increasing the two neighbouring entries by one:
\begin{align}
(c_1,...,c_i,c_{i+1},...,c_n) \rightarrow (c_1,...,c_{i}+1,1,c_{i+1}+1,...,c_n)
\label{eqn:inverse1}
\end{align}
(b) Inserting a -1 and decreasing the two neighbouring entries by one:
\begin{align}
(c_1,...,c_i,c_{i+1},...,c_n) \rightarrow (c_1,...,c_{i}-1,-1,c_{i+1}-1,...,c_n)
\label{eqn:inverse2}
\end{align}
(c) Splitting one entry and inserting a zero in between:
\begin{align}
(c_1,...,c_i,...,c_n) \rightarrow (c_1,...,c_{i}-z,0,z,...,c_n)
\label{eqn:inverse3}
\end{align}
or
\begin{align}
(c_1,...,c_i,...,c_n) \rightarrow (c_1,...,z,0,c_{i}-z,...,c_n)
\label{eqn:inverse4}
\end{align}
Note that in the third case, there are different inverses to the initial surgery operation since the expansion operation contains a parameter z. Since $c_i-z+z=c_i$ for any given $z$, the operation is indeed inverse to adding the two entries and removing the zero in between as in (\ref{eqn:surgery3}). Note also that strictly speaking, one of the equations (\ref{eqn:inverse3}) and (\ref{eqn:inverse4}) is redundant since replacing $z$ by $c_i-z$ in the lower one, one obtains again the upper one. However, it is easier to establish the relation between quiddity cycles and triangulations of polygons in the next chapter working with both equations.
\label{lm:inverse}
\end{bla}
\subsection{Continued Fraction Expressions and the Euclidean Algorithm}
\label{sec:neuf}

\begin{sz}
There exists a bijection between the set of rationals and the set of sequences of the form $(a_1,...,a_{2m})$ with $m \in \N$ where $a_1 \in \Z$ and $a_i \in \N$ for all $i \in \{2,...,n\}$.\\
A classical theorem (see for example section 2.1 in \cite{Andreescu}) states that for any rational $q$, there is a unique way to write $q$ in the form
\begin{equation}
q=a_1+\frac{1}{a_2+\frac{1}{a_3+\frac{1}{...+\frac{1}{a_n}}}}=:[a_1,...,a_n].
\end{equation}
More concretely, there are exactly two ways to express a rational as above, and exactly one of these expressions contains an even number of coefficients. This notation is called a continued fraction\index{continued fraction}.\\ Furthermore, it is easy to show that $q \in (a_1,a_1+1]$. Obviously, it is $q>a_1$. It is an easy exercise to show inductively starting with $a_n=a_{2m}$ that each denominator is larger than or equal to one, and hence $q \leq a_1+1$. Thus, $a_1=\lfloor q \rfloor$ can take any integer values while all other coefficients are only allowed to take positive values. The claimed bijection is
\begin{equation}
(a_1,...,a_{2m}) \mapsto [a_1,...,a_{2m}].
\end{equation}
\label{sz:cf}
\end{sz}

\begin{sz}
Any rational $q$ can be written uniquely in the form
\begin{equation}
q=d_1-\frac{1}{d_2-\frac{1}{d_3-\frac{1}{...-\frac{1}{d_k}}}}=:[[d_1,...,d_k]],
\label{eqn:cf}
\end{equation}
where $d_1 \in \Z$ and $d_i \geq 2$ for all $i \in \{2,...,n\}$, and if $[a_1,...,a_{2m}]=[[d_1,...,d_k]]$, the coefficients $d_i$ can be obtained via the equation
\begin{equation}
(d_1,...,d_k)=(a_1+1,\underbrace{2,...,2}_{a_2-1~ \textrm{times}},a_3+2,\underbrace{2,...,2}_{a_4-1~ \textrm{times}},...,a_{2m-1}+2,\underbrace{2,...,2}_{a_{2m}-1~ \textrm{times}}).
\label{eqn:Hirz}
\end{equation}
This form of writing a rational is called a negative continued fraction\index{negative continued fraction}.\\
This theorem is adopted from chapter 2 in \cite{Farey} as well as the following definition. For the origin of the formula (\ref{eqn:Hirz}), see \cite{Hirzebruch}, p.241 and \cite{Hirzebruch2}, eqns. (22) and (23). The uniqueness of the negative continued fraction follows by the uniqueness of the standard continued fraction.
\end{sz}

\begin{bla}{Lemma: Negative Euclidean Algorithm}
For two given coprime integers $a$ and $c \neq 0$, there exist unique integers $q_1,...,q_k$, $r_1,...,r_{k-1}$ with $q_2,...,q_k \geq 2$ and $0 < r_i <r_{i-1}$ for all $i \in \{1,...,k-1\}$ when setting $r_0=c$ such that the following equations are fulfilled.
\begin{equation}
\begin{aligned}
a&=q_1c-r_1\\
c&=q_2r_1-r_2\\
r_1&=q_3r_2-r_3\\
&\ldots\\
r_{k-3}&=q_{k-1}r_{k-2}-r_{k-1}\\
r_{k-2}&=q_kr_{k-1}-0.
\end{aligned}
\label{eqn:euclide}
\end{equation}
For the fraction of $a$ and $c$, it holds
\begin{equation}
\frac{a}c=[[q_1,...,q_k]].
\label{eqn:ncfe}
\end{equation}
\end{bla}
\begin{bw}
Divide the first equation of (\ref{eqn:euclide}) by $c$ and obtain
\begin{equation}
\frac{a}c=q_1-\frac{r_1}c=q_1-\frac{1}{c/r_1},
\end{equation}
where the last step is done to obtain a form where the second equation of (\ref{eqn:euclide}) can be inserted the same way after being divided by $r_1$:
\begin{equation}
\frac{a}c=q_1-\frac{1}{c/r_1},=q_1-\frac{1}{q_2-\frac{r_2}{r_1}}=q_1-\frac{1}{q_2-\frac{1}{r_1/r_2}}.
\end{equation}
Iterating this process making successively use of all the equations of (\ref{eqn:euclide}), one finally obtains
\begin{equation}
\frac{a}c=[[q_1,...,q_k]].
\end{equation}
The conditions $0 < r_{i} < r_{i-1}$ for all $i \in \{1,...,k-1\}$ allow to deduce successively $q_1=\left \lceil \frac{a}c \right \rceil$, $q_2=\left \lceil \frac{c}{r_1} \right \rceil$ and $q_i=\left \lceil \frac{r_{i-2}}{r_{i-1}} \right \rceil$ for $i > 2$. Hence, the coefficients $q_i$ are uniquely defined, and since $0 < r_{i} < r_{i-1}$, it is $q_i \geq 2$ for all $i \geq 2$. If $r_i=0$, it follows $i=k$, $q_k=\left \lceil \frac{r_{k-2}}{r_{k-1}} \right \rceil=\frac{r_{k-2}}{r_{k-1}}$, and the algorithm comes to an end. The uniqueness of the coefficients $q_i$ also induces uniqueness for the coefficients $r_i$.
\end{bw}
\section{Expressing $PSL(2,\Z)$ in Terms of Continuant Matrices}
\label{sec:un}

In chapter 7 of \cite{Farey}, Morier-Génoud and Ovsienko developped a method to obtain a presentation for matrices in $PSL(2,\Z)$ in terms of continuant matrices with only positive integers as coefficients and imposed conditions for this presentation to be minimal in terms of its length. Concretely, they found out that a matrix $M(1)$ can only occur at the first two or at the last two spots in the presentation.\\
Now, the goal is to investigate what happens if negative (and zero) coefficients are also allowed. It turns out that all entries with absolute value smaller than three play a crucial role for this purpose.\\
Developping this general minimal presentation of $PSL(2,\Z)$ in terms of continuant matrices with integer coefficients needs some preparation. It turns out to be useful for a given matrix $A \in PSL(2,\Z)$ to take a small detour via a presentation with continuant matrices that is longer than the minimal one in general, but obtainable with the help of the negative continued fraction expansion of the ratio of the two left entries of the matrix. The process to develop this presentation is described in this section.\\
For simplification of notation, in the rest of this section, a matrix $A$ stands for its equivalence class $[A]$ in $PSL(2,\Z)$ unless indicated otherwise.\\
The first subsection is dedicated to the formulation and proof of theorem \ref{sz:calc2} which is the central theorem of this section. An idea developed in \cite{Conrad} to find a well-defined presentation for each $A \in PSL(2,\Z)$ similar to the presentation given in chapter 7 of \cite{Farey} is used for the proof. In the second subsection, some examples for and consequences of this theorem are given.

\subsection{Derivation of a Presentation Using Negative Continued Fractions}

\begin{lm}
Let $A=\begin{pmatrix} a&b \\ c&d \end{pmatrix} \in SL(2,\Z)$ with $|a|>1$, $c \neq 0$, and let $a$ and $c$ have the same sign. Then $\left \lceil \frac{b}a \right \rceil = \left \lceil \frac{d}c \right \rceil$.
\label{lm:determinant}
\end{lm}
\begin{bw}
Assume $\left \lceil \frac{b}a \right \rceil \neq \left \lceil \frac{d}c \right \rceil$. It is $\frac{b}a < \frac{d}c$ since otherwise det$(A)$ would be negative, and hence follows $\left \lceil \frac{b}a \right \rceil < \left \lceil \frac{d}c \right \rceil$. Then there is an integer $n$ such that $\frac{b}a \leq n < \frac{d}c$. Multiplying by $ac$ yields
\begin{equation}
bc \leq nac < ad.
\label{eqn:determinant}
\end{equation}
All three terms in equation (\ref{eqn:determinant}) are integers. Furthermore $1=$det$(A)=ad-bc$, and the right inequality is strict. Therefore, $bc=nac$ resp. $b=na$. So, $a|b$, and hence, $a=$gcd$(a,b)$. But gcd$(a,b)|$det$(A)$, and therefore, $a=\pm 1$, a contradiction.
\end{bw}

\begin{sz}
Let $A=\begin{pmatrix} a&b \\ c&d \end{pmatrix} \in PSL(2,\Z)$. If $c \neq 0$, it is
\begin{equation}
A=M(q_1,...,q_k,m+1,1,1)=M(q_1,...,q_k,m,0)
\label{eqn:presente}
\end{equation}
with
\begin{equation}
\frac{a}c=[[q_1,...,q_k]]~~\textrm{and}~~m=\left \lceil \frac{d}c \right \rceil.
\end{equation}
 If $c=0$, it is $A=M(m+1,1,1)=M(m,0)$ with $m=\frac{b}a=\frac{b}d$.\\
The presentation given in equation (\ref{eqn:presente}) is unique as a presentation of $A$ ending on zero and fulfilling the condition $q_2,...,q_k \geq 2$, i.e. that all entries except the first one and the two last ones have a value of at least two.\\ \\
Regarding the first expression in equation (\ref{eqn:presente}), note that $q_1$ and $m+1$ are the only integers in the sequence that can take all integers as their value. Especially, $q_1 \in \{0,1\}$ if $|a|<|c|$ or $a=c=1$ and $q_1 \notin \{0,1\}$ otherwise, and $q_1 \geq 2$ if and only if $a$ and $c$ have the same sign and $|a|>|c|$.\\
Note also that in the first expression, all coefficients are positive if and only if $\left \lceil \frac{d}c \right \rceil$ and $\left \lceil \frac{a}c \right \rceil$ are positive, i.e. $a$, $c$ and $d$ (and hence also $b$ if non-zero) all have the same sign, or, in the case $c=0$, $a$, $b$ and $d$ have the same sign. In other words, all coefficients are positive if and only if all (non-zero) entries of $A$ have the same sign.
\label{sz:calc2}
\end{sz}
\begin{bw}
(I) \textit{Existence}: Let $c \neq 0$ and $S$ and $T$ denote
\begin{equation}
S=\begin{pmatrix} 0&-1 \\ 1&0 \end{pmatrix},~~T=\begin{pmatrix} 1&1 \\ 0&1 \end{pmatrix},
\end{equation}
and hence, for $m \in \Z$,
\begin{equation}
-M(m+1,1,1)=-M(m,0)=T^m=\begin{pmatrix} 1&m \\ 0&1 \end{pmatrix}
\label{eqn:wurst}
\end{equation}
and $M(m)=T^mS$. Since $S^{-1}=-S$, it is also $M(m)^{-1}=-ST^{-m}$. Since $A \in PSL(2,\Z)$, assume $c \geq 0$, otherwise consider $-A$. The goal is now to prove equation (\ref{eqn:presente}) by showing
\begin{equation}
(M(q_1,...,q_k,m,0))^{-1}A=Id.
\label{eqn:importante}
\end{equation}
To do so, set $A_0=A$ and
\begin{equation}
A_{i+1}=\begin{pmatrix} a_{i+1}&b_{i+1}\\c_{i+1}&d_{i+1} \end{pmatrix}:=M(q_{i+1})^{-1}A_i=-ST\overset{-q_{i+1}}{}A_i=\begin{pmatrix} c_i&d_i\\q_{i+1}c_i-a_i&q_{i+1}d_i-b_i \end{pmatrix}
\end{equation}
for $i \in \{0,...,k\}$.
Now, one investigates what happens when the matrix $A$ is subsequently multiplied on the left with these blocks of the form $(M(q_i))^{-1}=-ST\overset{-q_i}{}$. It is
\begin{equation}
A_1=-ST\overset{-q_1}{}A=\begin{pmatrix} c&d \\ q_1c-a&q_1d-b \end{pmatrix}=\begin{pmatrix} c&d \\ r_1&q_1d-b \end{pmatrix}
\end{equation}
Applying $ST\overset{-q_2}{}$ from the left then yields
\begin{equation}
A_2=(-ST\overset{-q_2}{})(-ST\overset{-q_1}{})A=\begin{pmatrix} r_1&\ast \\ q_2r_1-c&\ast \end{pmatrix} = \begin{pmatrix} r_1&\ast \\ r_2&\ast \end{pmatrix}
\end{equation}
where the second column of the matrix is now given by $\begin{pmatrix} \ast \\ \ast \end{pmatrix} = (-ST\overset{-q_2}{})(-ST\overset{-q_1}{})\begin{pmatrix} b \\ d \end{pmatrix}$. Iterating this procedure yields
\begin{equation}
A_i=(-ST\overset{-q_i}{})...(-ST\overset{-q_1}{})A=\begin{pmatrix} r_{i-1}&\ast \\ r_i&\ast \end{pmatrix}
\end{equation}
and finally
\begin{equation}
A_k=(-ST\overset{-q_k}{})...(-ST\overset{-q_1}{})A=\begin{pmatrix} r_{k-1}&\ast \\ 0&\ast \end{pmatrix}.
\label{eqn:longmulti}
\end{equation}
Now, $\mathrm{det}(A_k)=1$ and $\mathrm{gcd}(a_k,c_k)|\mathrm{det}(A_k)$ implies $\mathrm{gcd}(a_k,c_k)=1$, and hence, $r_{k-1}=1$. Since the lower left entry of the matrix in equation (\ref{eqn:longmulti}) is zero and the upper left entry is 1, the lower left entry must also be 1 to ensure that the determinant is still 1. Conclude that for some $m \in \Z$, it is
\begin{equation}
A_k=(-ST\overset{-q_k}{})...(-ST\overset{-q_1}{})A=\begin{pmatrix} r_{k-1}&\ast \\ 0&\ast \end{pmatrix}=\begin{pmatrix} 1&m \\ 0&1 \end{pmatrix} = T^m.
\label{eqn:minusst}
\end{equation}
Thus, applying $T^{-m}=(M(m,0))^{-1}$ (in $PSL(2,\Z)$, the minus sign in equation (\ref{eqn:wurst}) can be skipped) from the left yields the identity matrix, equation (\ref{eqn:importante}) is indeed fulfilled, and equation (\ref{eqn:presente}) holds for some $m \in \Z$. It remains to show that $m=\left \lceil \frac{d}c \right \rceil$.\\ \\
In the case $k=1$, it is
\begin{equation}
A_k=A_1=\begin{pmatrix} c&d \\ r_1&q_1d-b \end{pmatrix}=\begin{pmatrix} 1&d \\ 0&1 \end{pmatrix}.
\end{equation}
It follows $m=d=\frac{d}1=\frac{d}c=\left \lceil \frac{d}c \right \rceil$.\\
In the case $k \geq 2$, by construction, for $1 \leq i \leq k-1 = \max\{i|c_i>0\}$, it is $a_i=r_{i-1}>r_i=c_i > 0$. Hence, for $1 \leq i \leq k-1$, the matrix $A_i$ fulfils the conditions of lemma \ref{lm:determinant}, and it is
\begin{equation}
\left \lceil \frac{b_{i+1}}{a_{i+1}} \right \rceil = \left \lceil \frac{d_i}{c_i} \right \rceil = \left \lceil \frac{b_i}{a_i} \right \rceil.
\end{equation}
By iteration, it follows
\begin{equation}
m = \left \lceil \frac{m}{1} \right \rceil = \left \lceil \frac{b_k}{a_k} \right \rceil = \left \lceil \frac{d_{k-1}}{c_{k-1}} \right \rceil = \left \lceil \frac{b_{k-1}}{a_{k-1}} \right \rceil = ... = \left \lceil \frac{d_1}{c_1} \right \rceil = \left \lceil \frac{b_1}{a_1} \right \rceil = \left \lceil \frac{d}{c} \right \rceil
\end{equation}
since, recalling equation (\ref{eqn:minusst}), it is $A_k=(-ST\overset{-q_k}{})...(-ST\overset{-q_1}{})A=\begin{pmatrix} 1&m\\0&1 \end{pmatrix}$.\\ \\
Finally, if $c=0$, $\det A=1$ implies $A=\pm \begin{pmatrix} 1&m\\0&1 \end{pmatrix}=M(m,0)$ in $PSL(2,\Z)$ for some $m \in \Z$, and it follows $\frac{b}a=\frac{b}d=\frac{m}1=m$.\\ \\
(II) \textit{Uniqueness}: Suppose $A=M(q_1,...,q_k,m,0)$ with $q_2,...,q_k \geq 2$ and set\\
$A_i=\begin{pmatrix} a_i&b_i \\ c_i&d_i \end{pmatrix} := \prod_{s=i}^k M(q_s) =  M(q_i,...,q_k)$ for $i \in \{1,...,k\}$. Then, it is
\begin{equation}
A_k=\begin{pmatrix} q_k&-1 \\ 1&0 \end{pmatrix}
\end{equation}
and
\begin{equation}
A_{i-1}=\begin{pmatrix} q_{i-1}&-1 \\ 1&0 \end{pmatrix} \begin{pmatrix} a_i&b_i \\ c_i&d_i \end{pmatrix} = \begin{pmatrix} q_{i-1}a_i-c_i&q_{i-1}b_i-d_i \\ a_i&b_i \end{pmatrix} = \begin{pmatrix}a_{i-1}&b_{i-1}\\c_{i-1}&d_{i-1} \end{pmatrix}~\textrm{for}~i \in \{2,...,k\}.
\end{equation}
If $k=1$, then $A_k=A_1$, $c_1=1 \neq 0$, and $\frac{a}c=a=q_k=[[q_1]]$ follows. Otherwise, $k \geq 2$ and $a_k=q_k > c_k=1>0$. Then, it follows $[[q_{k-1},q_k]]=q_{k-1}-\frac{1}{q_k}=\frac{q_{k-1}q_k-1}{q_k}=\frac{a_{k-1}}{c_{k-1}}$, $c_{k-1}=a_k>0$, and if $k>2$, also $a_{k-1}=q_{k-1}a_k-c_k>c_{k-1}=a_k$ because it is $q_{k-1} \geq 2$ and $c_k < a_k$.\\
Now, iterate and deduce that for all $i > 2$, it follows inductively
\begin{equation}
[[q_{i-1},...,q_k]]=q_{i-1}-\frac{1}{[[q_i,...,q_k]]}=q_{i-1}-\frac{c_i}{a_i}=\frac{q_{i-1}a_i-c_i}{a_i}=\frac{a_{i-1}}{c_{i-1}},
\end{equation}
and $a_{i-1}>c_{i-1}>0$. Hence, especially $\frac{a_2}{c_2}=[[q_2,...,q_k]]$ and $a_2>c_2>0$. Then, it follows also
\begin{equation}
[[q_1,...,q_k]]=q_1-\frac{1}{[[q_2,...,q_k]]}=q_1-\frac{c_2}{a_2}=\frac{q_1a_2-c_2}{a_2}=\frac{a_1}{c_1}
\end{equation}
and $c_1 = a_2 > 0$, hence especially $c_1 \neq 0$.\\ \\
Now, consider the matrix $A_1=M(q_1,...,q_k)=\begin{pmatrix} a_1&b_1 \\ c_1&d_1 \end{pmatrix}$ and multiply $M(m,0)=\begin{pmatrix} 1&m \\ 0&1 \end{pmatrix}$ (in $PSL(2,\Z)$) from the right to obtain $A$. Since this multiplication does only affect the right column of the matrix, it is $a=a_1$, $c=c_1$, and hence follows $c \neq 0$ and $\frac{a}c=\frac{a_1}{c_1}=[[q_1,...,q_k]]$.\\
So, there is only one possible choice for each of the values $q_1,...,q_k$, i.e. they are uniquely determined. Consider now the matrices $M(q_1,...,q_k,m,0)$ and $M(q_1,...,q_k,m',0)$. Since all of the matrices in the two products except $M(m)$ resp. $M(m')$ are the same, it follows also $M(m)=M(m')$, and hence $m=m'$. As shown above, $M(q_1,...,q_k,\left \lceil \frac{d}c \right \rceil,0)=A$ is a presentation of $A$, and hence this is the only possibility to fulfil the conditions above, i.e. it follows $m=\left \lceil \frac{d}c \right \rceil$.
\end{bw}

\subsection{Examples and Corollaries}

\begin{bsp}
Consider the matrix $A=\begin{pmatrix} 17&12\\7&5 \end{pmatrix}$. Then, with $q_1=\left \lceil \frac{17}7 \right \rceil=3$, and then by iteration $q_2=2$ and $q_3=4$, it is
\begin{equation}
\begin{aligned}
A_1&:=-ST^{-3}A= \begin{pmatrix} 7&5\\4&3 \end{pmatrix}\\
A_2&:=-ST^{-2}A_1= \begin{pmatrix} 4&3\\1&1 \end{pmatrix}\\
A_3&:=-ST^{-4}A_2= \begin{pmatrix} 1&1\\0&1 \end{pmatrix}.
\end{aligned}
\end{equation}
Hence, $m=1$ and $\frac{17}7=[[3,2,4]]$, and therefore $A=M(3,2,4,1,0)=M(3,2,3,-1)$ in $PSL(2,\Z)$ where the last equality follows by equation (\ref{eqn:calc1}).
\end{bsp}

\begin{bsp}
(a) Consider the matrix $A=\begin{pmatrix} -60&13\\23&-5 \end{pmatrix}$. It is $\frac{-60}{23}=[[-2,2,3,5]]$ and $\left \lceil \frac{-5}{23} \right \rceil=0$. Hence, $A=M(-2,2,3,5,0,0)=M(-2,2,3,5)$ in $PSL(2,\Z)$.\\
(b) Consider the matrix $B=\begin{pmatrix} -4&-9\\9&20 \end{pmatrix}$. It is $\frac{-4}{9}=[[0,3,2,2,2]]$ and $\left \lceil \frac{20}{9} \right \rceil=3$. Hence, $B=M(0,3,2,2,2,3,0)$ in $PSL(2,\Z)$.
\end{bsp}

\begin{rem}
In the case $c=0$, one can also interpret $a/c$ as infinitely large and interpret the empty negative continued fraction $[[]]$ as an expression for it. Then, with $k=0$, the length of $[[]]$, equation (\ref{eqn:presente}) is still valid.
\end{rem}

\vspace{5mm}

Now, it is possible to deduce the following results which are also derived in \cite{Farey} (thm. 7.3 (ii) and cor. 7.5). Part (iii) is not a direct consequence of theorem \ref{sz:calc2}, but follows from part (ii). The proof for part (iii) which is included in this text for the sake of completeness is essentially taken from \cite{Farey}.

\begin{cor}
\begin{itemize}
\item (i) Let $A=\begin{pmatrix} a&b \\ c&d \end{pmatrix} \in PSL(2,\Z)$, $a,b,c,d>0$ and $a>b$. Then $A=M(q_1,...,q_k,1,0)=M(q_1,...,q_k,2,1,1)$ with $[[q_1,...,q_k]]=a/c$.
\item (ii) Let $A=\begin{pmatrix} a&-b \\ c&-d \end{pmatrix} \in PSL(2,\Z)$ with $a,b,c,d>0$ and $a>b$. Then $A=M(q_1,...,q_k)$ with $\frac{a}c=[[q_1,...,q_k]]$.
\item (iii) Let $A=\begin{pmatrix} a&b \\ c&d \end{pmatrix} \in PSL(2,\Z)$ with $a,b,c,d>0$ and $a<b$. Then\\
$A=M(p_1,...,p_k,0)=M(p_1,...,p_{k-1},p_k+1,1,1)$ with $\frac{b}d=[[p_1,...,p_k]]$.
\end{itemize}
\label{cor:misone}
\end{cor}
\begin{bw}
\begin{itemize}
\item (i) Since $a>b>0$, it follows $a>1$, and hence lemma \ref{lm:determinant} is applicable. Thus, one obtains $\left \lceil \frac{d}c \right \rceil = \left \lceil \frac{b}a \right \rceil=1$ where the last equality is true because $0<b<a$. Now, the statement follows with theorem \ref{sz:calc2}.
\item (ii) Since $a>b$, it must be $c>d$ to enable still det$(A)=bc-ad=1$. Hence, $m=\left \lceil \frac{-d}c \right \rceil=0$, and by theorem \ref{sz:calc2}
\begin{equation}
A=M(q_1,...,q_k,0,0)=M(q_1,...,q_k)M(0,0)=M(q_1,...,q_k)
\end{equation}
in $PSL(2,\Z)$ since $(0,0)$ is a quiddity cycle what is equivalent to $M(0,0)=\pm Id$.
\item (iii) For $A$ as in the conditions, the matrix
\begin{equation}
AS=\begin{pmatrix} a&b \\ c&d \end{pmatrix}\begin{pmatrix} 0&-1 \\ 1&0 \end{pmatrix}=\begin{pmatrix} b&-a \\ d&-c \end{pmatrix}
\end{equation}
fulfils the conditions of part (ii). It follows
\begin{equation}
\begin{aligned}
A=(AS)S^{-1}&=M(p_1,...,p_k) (-M(0))\\
&=M(p_1,...,p_k,0)\\
&=M(p_1,...,p_{k-1},p_k+1,1,1)
\end{aligned}
\end{equation}
in $PSL(2,\Z)$.
\end{itemize}
\end{bw}
\section{Uniqueness and Minimality of the Presentation}
\label{sec:deux}
In the last section, an algorithm to find a presentation with continuant matrices $M(c_1,...,\\c_{k+3})=M(q_1,...,q_k,m+1,1,1)=M(q_1,...,q_k,m,0)$ for any matrix $A \in PSL(2,\Z)$ was derived. In the first expression, all of the coefficients in the sequence are positive except maybe the first one and $m+1$. However, allowing also negative coefficients at all positions, one may find a shorter presentation as example \ref{ex:notshort} shows. The goal of this section is now to create an algorithm to find the shortest possible presentation of an arbitrary matrix $A \in PSL(2,\Z)$, i.e. a presentation $A=M(c_1,...,c_k)$ with $k$ as small as possible. Furthermore, it is examined under which circumstances this minimal presentation is unique.\\
In the first subsection, it is shown that presentations including consecutive twos or -2's in their inner part cannot be minimal. In the second subsection, it becomes clear that any presentation of a matrix $A \in PSL(2,\Z)$ can be transformed into the unique one of section \ref{sec:un} using only the transformations introduced previously. In the third subsection, it is dealt with subsequences of the form $(\pm 2,\pm 3,...,\pm 3,\pm 2)$. In the fourth subsection, finally the criteria for a presentation of a matrix $A \in PSL(2,\Z)$ to be minimal and for such a minimal presentation to be unique are stated and proven. It will turn out that the presentation becomes unique with the canonical choice that except perhaps as first or last coefficient, $c_i=-2$ is disallowed.\\
In the whole section, eventual minus signs in front of matrices are skipped since the considered group is $PSL(2,\Z)$.

\subsection{Subsequences of Consecutive Twos or -2's}

\begin{bsp}
Consider the matrix
\begin{equation*}
A=\begin{pmatrix} 38&-17 \\ 9&-4 \end{pmatrix} = M(5,2,2,2,3)= M(5,2,2,2,3,0,0) \in PSL(2,\Z).
\end{equation*}
However, calculations show that also $M(4,-4,2)=A$.
\label{ex:notshort}
\end{bsp}

\begin{lm}
For $k \in \N_0$, it is
\begin{equation}
(M(2))^k=\begin{pmatrix} k+1&-k \\ k&-(k-1) \end{pmatrix}~~\textrm{and}~~(M(-2))^k= (-1)^k\begin{pmatrix} k+1&k \\ -k&-(k-1) \end{pmatrix}.
\end{equation}
In $PSL(2,\Z)$, the factor $(-1)^k$ in the second equation can also be skipped.
\end{lm}
\begin{bw}
Induction over $k$. Clearly, for $k=0$, the statement is true. Assume now that for $k-1$, it is already shown
\begin{equation}
(M(2))^{k-1}=\begin{pmatrix} k&-(k-1) \\ k-1&-(k-2) \end{pmatrix}~~\textrm{and}~~(M(-2))^{k-1}= (-1)^{k-1}\begin{pmatrix} k&k-1 \\ -(k-1)&-(k-2) \end{pmatrix}.
\end{equation}
Then, it follows
\begin{equation}
(M(2))^k= \begin{pmatrix} 2&-1 \\ 1&0 \end{pmatrix} \begin{pmatrix} k&-(k-1) \\ k-1&-(k-2) \end{pmatrix} = \begin{pmatrix} k+1&-k \\ k&-(k-1) \end{pmatrix}
\end{equation}
and
\begin{equation}
(M(-2))^k= (-1)^{k-1} \begin{pmatrix} -2&-1 \\ 1&0 \end{pmatrix} \begin{pmatrix} k&k-1 \\ -(k-1)&-(k-2) \end{pmatrix} = (-1)^k \begin{pmatrix} k+1&k \\ -k&-(k-1) \end{pmatrix}.
\end{equation}
\end{bw}

\begin{cor}
For $k \in \N$, the following identities hold:
\begin{equation}
M(a) (M(2))^k M(b)= -M(a-1,-(k+1),b-1)
\label{eqn:powerone}
\end{equation}
\begin{equation}
M(a)(M(-2))^k M(b)=(-1)^k M(a+1,k+1,b+1)
\label{eqn:powertwo}
\end{equation}
\begin{equation}
M(a,-k)=-M(a,-k,0,0)=M(a+1) M(2)^{k-1} M(1,0).
\label{eqn:powerthree}
\end{equation}
Note again that in $PSL(2,\Z)$, the signs can also be skipped, and that the first two equations also hold for $k=0$. In this special case, they turn precisely into the equations (\ref{eqn:calc2}) resp. (\ref{eqn:calc1}).
\label{cor:powertwo}
\end{cor}
\begin{bw}
Using the previous lemma, this is an easy calculation.
\end{bw}

\begin{bla}{Another Type of Surgery Operations}
Corollary \ref{cor:powertwo} allows to introduce a new type of surgery operations\index{surgery operations}. Considering a Matrix $A=M(c_1,...,c_k) \in PSL(2,\Z)$ where the sequence contains a coefficient $-k:=c_i \leq -1$, the transformation
\begin{align}
c:= (c_1,...,c_{i-1},-k,c_{i+1},...,c_n) \rightarrow (c_1,...,c_{i-1}+1,\underbrace{2,...,2}_{(k-1)~\textrm{times}},c_{i+1}+1,...,c_n) =:\tilde c
\label{eqn:surgerytwo1}
\end{align}
can be applied, and by equation (\ref{eqn:powerone}), it is also $A=M(c)=M(\tilde c)$. For $k \geq 3$, this operation increases the number of entries of the sequence, and hence, its inverse operation can be seen as another surgery operation. For $k=2$, the number of entries of the sequence remains the same, and for $k=1$, the operation is nothing else than operation (\ref{eqn:surgery2}).\\
If $c$ contains a subsequence $(c_{i+1},...,c_{i+k})=(-2,...,-2)$ with $k \geq 0$, use equation (\ref{eqn:powertwo}) to deduce that it is also $A=M(c)=M(\tilde c)$ when transforming
\begin{align}
c:= (c_1,...,c_i,\underbrace{-2,...,-2}_{k~\mathrm{times}},c_{i+k+1},...,c_n) \rightarrow (c_1,...,c_i+1,k+1,c_{i+k+1}+1,...,c_n) =:\tilde c.
\label{eqn:surgerytwo2}
\end{align}
For $k=0$, this is precisely the inverse operation to (\ref{eqn:surgery1}). For $k=1$, the number of entries stays unchanged, and for $k \geq 2$, it is reduced. Hence, in this case, the transformation can also be viewed as a surgery operation. Furthermore, if $-k:=c_n \leq -1$, the transformation
\begin{align}
c:= (c_1,...,c_{n-1},-k) \rightarrow (c_1,...c_{n-1}+1,\underbrace{2,...,2}_{(k-1)~\textrm{times}},1,0) =:\tilde c
\label{eqn:surgerytwo3}
\end{align}
can be applied, and by equation (\ref{eqn:powerthree}), $A=M(c)=M(\tilde c)$ still holds. The inverse operation is also called a surgery operation since it reduces the number of entries of the sequence as well.
\label{bla:newsurgery}
\end{bla}

\subsection{Relation to the Unique Presentation from Section \ref{sec:un}}

\begin{sz}
Let $A=M(c_1,...,c_n) \in PSL(2,\Z)$. Then, with a finite number of the transformations described in lemma \ref{lm:surg} and \ref{bla:newsurgery} and - if necessary - once adding $(0,0)$ at the end, the sequence $(c_1,...,c_k)$ can be transformed into a sequence $(q_1,...,q_k,m,0)$ such that $q_i \geq 2$ for all $i \in \{2,...,k\}$ and $A=M(q_1,...,q_k,m,0)$.
\label{sz:sequence}
\end{sz}
\begin{bw}
1) First, eliminate all entries $c_i$ with $|c_i| \leq 1$ that are not at the beginning or end of the sequence by applying transformations from lemma \ref{lm:surg}. Since $n$ is finite and all these transformations reduce the length of the sequence, after a finite number of steps, one obtains a sequence where all entries $c_i$ except perhaps the first and the last one fulfil $|c_i| \geq 2$.\\
2) The second step is eliminating all entries that are -2 except perhaps the first and the last entry of the sequence. To do this, apply transformations of type (\ref{eqn:surgerytwo2}) with suitable values of $k$. The transformations either do not change the length of the sequence (for $k=1$) or reduce it (for $k>1$), and each of the transformation reduces the number of negative entries in the sequence by $k$.\\
It is also impossible that such a transformation creates a new entry $c_i$ with $|c_i|<2$ in the sequence except perhaps at the beginning or at the end if $k$ is chosen maximal since if one of the neighbouring entries of the subsequence $(-2,...,-2)$ is also $-2$, this entry can also be included in the subsequence, the value for $k$ can be increased by one, and there will appear no new $-1$ by increasing the neighbor entry by one. So, after a finite number of operations, for all entries $c_i$ except perhaps the first and the last one, it is either $c_i \leq -3$ or $c_i > 2$.\\
3) Now, replace all entries $c_i$ with $c_i \leq -3$ except perhaps the first and the last entry of the sequence by applying a transformation of type (\ref{eqn:surgerytwo1}) with $k=|c_i|$ with a subsequence of $(k-1)$ twos. Since there has to be a 1 added to the two neighbor entries of $c_i$, if $c_{i-1}=-3$ or $c_{i+1}=-3$, applying this transformation produces a new entry $-2$, and one has to go back to step 2. However, each transformation of type (\ref{eqn:surgerytwo1}) also reduces the total number of negative entries by one, and hence, after a finite number of steps, all entries except perhaps the first and the last one are at least 2.\\
4) If the last entry is zero, the desired sequence $(q_1,...,q_k,m,0)$ is already found. If the last entry is negative, apply transformation (\ref{eqn:surgerytwo3}) to obtain a sequence $(q_1,...,q_k,m,0)$ fulfilling all conditions above. If the last entry is larger than one, append $(0,0)$ to the sequence (which is possible since $M(0,0)=-Id=Id$ in $PSL(2,\Z)$), and we are also done.\\
The most complicated case is when the last entry is one. Appending $(0,0)$ is not possible since otherwise the entry that becomes the third last entry is one what is not allowed.\\
But it is possible to apply equation (\ref{eqn:powerone}) with $b=1$, and then the 1 at the end is decreased by one and is hence zero. If the sequence $(c_1,...,c_r,2,...,2,1)$ is not of the form $(2,2,...,2,2,1)$, i.e. does not only consist of twos and the one at the end, choose for $k$ the number of subsequent twos directly before the one at the end, and if there are none, $k=0$. After transformation, the sequence turns into $(c_1,...,c_r-1,-(k+1),0)$. Since $c_r$ was the last entry larger than two of the initial sequence, it is $c_r-1 \geq 2$, and hence, the resulting sequence $(c_1,...,c_r-1,-(k+1),0)$ is the desired sequence $(q_1,...,q_k,m,0)$.\\
 If the sequence is of the form $(2,2,...,2,2,1)$ and contains $l$ twos, choose $k=l-1$, and the first 2 is the $a$ in equation (\ref{eqn:powerone}), and one obtains the sequence $(1,m,0)$ with $m=-l$.\\ \\
Since all the used transformations do not change the represented matrix\index{represented matrix} (modulo a global sign), it follows $A=M(q_1,...,q_k,m,0)$.
\end{bw}

\begin{cor}
Let $A=\begin{pmatrix} a&b \\ c&d \end{pmatrix} \in PSL(2,\Z)$. Let furthermore $\mathcal M$ be the set of all operations of the types (\ref{eqn:surgery1})-(\ref{eqn:surgery3}) and (\ref{eqn:surgerytwo1})-(\ref{eqn:surgerytwo3}) and all operations that are inverse to one of those operations. Note that the operation to append $(0,0)$ at the end of the cycle does not need to be included here separately since this can also be viewed as the inverse operation to transformation (\ref{eqn:surgery3}) with $i=n-1$ and $c_n=0$. Also, the transformation (\ref{eqn:surgerytwo3}) is a special case of an inverse of the transformation (\ref{eqn:surgery1}) executed after appending $(0,0)$.\\
Then, for $c \neq 0$, there is a unique presentation of $A$ of the form $M(q_1,...,q_k,m,0)$ with $q_2,...,q_k \geq 2$, and this presentation can be found by calculating $(q_1,...,q_k)$ via the equations $[[q_1,...,q_k]]=\frac{a}c$ and $m=\left \lceil \frac{d}c \right \rceil$, and, if $c=0$, there is a unique presentation of $A$ which is $A=M(m,0)$ with $m=\left \lceil \frac{b}a \right \rceil=\left \lceil \frac{b}d \right \rceil$.\\
Furthermore, if $M(c_1,...,c_k)=M(d_1,...,d_l)$ are two different presentations of $A$, it is possible to transform the sequence $(c_1,...,c_k)$ into the sequence $(d_1,...,d_l)$ with a finite number of transformations out of $\mathcal M$.
\label{cor:unice}
\end{cor}
\begin{bw}
By theorem \ref{sz:calc2}, it follows that $M(q_1,...,q_k,m,0)$ as above is indeed a presentation of $A$. For $c \neq 0$, the uniqueness of this presentation follows by theorem \ref{sz:calc2}, and for $c=0$, it is trivial since for $m \neq m'$, obviously it follows $M(m,0) \neq M(m',0)$.\\
For the second part of the statement, first transform the sequence $(c_1,...,c_k)$ into the unique sequence of the form $(q_1,...,q_k,m,0)$ with $q_2,...,q_k \geq 2$ that represents the same matrix as described in the proof of theorem \ref{sz:sequence}, and then continue by transforming $(q_1,...,q_k,m,0)$ into $(d_1,...,d_l)$ by applying the inverse procedure.
\end{bw}

\subsection{Subsequences of the Type $(\pm 2,\pm 3,...,\pm 3,\pm 2)$}

\begin{lm}
It is
\begin{equation}
\begin{aligned}
&M(a,2,\underbrace{3,...,3}_{l~\mathrm{times}},2,b)=M(a-1,\underbrace{-3,...,-3}_{(l+1)~\mathrm{times}},b-1)~\textrm{and}\\
&M(a,-2,\underbrace{-3,...,-3}_{l~\mathrm{times}},-2,b)=M(a+1,\underbrace{3,...,3}_{(l+1)~\mathrm{times}},b+1)
\end{aligned}
\end{equation}
or explicitly
\begin{equation}
\begin{aligned}
&M(a,2)M(3)^lM(2,b)=M(a-1)M(-3)^{l+1}M(b-1)~\textrm{and}\\
&M(a,-2)M(-3)^lM(-2,b)=M(a+1)M(3)^{l+1}M(b+1),
\end{aligned}
\end{equation}
and hence, there are two more types of possible surgery operations\index{surgery operations}, namely
\begin{equation}
\begin{aligned}
&(a,2,\underbrace{3,...,3}_{l~\mathrm{times}},2,b) \to (a-1,\underbrace{-3,...,-3}_{(l+1)~\mathrm{times}},b-1)~\textrm{and}\\
&(a,-2,\underbrace{-3,...,-3}_{l~\mathrm{times}},-2,b) \to (a+1,\underbrace{3,...,3}_{(l+1)~\mathrm{times}},b+1),
\end{aligned}
\label{eqn:232}
\end{equation}
letting the represented matrix of the sequence invariant and reducing the total number of entries by one.
\label{lm:232}
\end{lm}
\begin{bw}
To prove the first equation, apply iteratively equation (\ref{eqn:powerone}) $l$ times with $k=1$ and finally once with $k=2$ as follows:
\begin{equation}
\begin{aligned}
&M(a)M(2)M(3)^lM(2)M(b)\\
&=M(a-1)M(-2)M(2)M(3)^{l-1}M(2)M(b)\\
&=M(a-1)M(-3)M(-2)M(2)M(3)^{l-2}M(2)M(b)\\
&=M(a-1)M(-3)^2M(-2)M(2)M(3)^{l-3}M(2)M(b)\\
&=...=M(a-1)M(-3)^{l-1}M(-2)M(2)M(2)M(b)\\
&=M(a-1)M(-3)^{l+1}M(b-1).
\end{aligned}
\end{equation}
The second equation follows analogously when applying equation (\ref{eqn:powertwo}) instead of equation (\ref{eqn:powerone}).
\end{bw}

\subsection{Criteria for Minimality and Uniqueness of the Presentation and Proof of the Statement}

\begin{bla}{Minimality Criteria\index{minimality criteria}}
Now, it is possible to state some criteria for a given presentation $A=M(c_1,...,c_n)$ to conclude whether it is possible to shorten it or not, i.e. whether it is minimal or not.\\
\begin{itemize}
\item If the subsequence $(c_2,...,c_{n-1})$ contains entries with absolute value smaller than two, the sequence cannot be minimal since then the number of entries can be reduced by applying one of the transformations (\ref{eqn:surgery1})-(\ref{eqn:surgery3}).
\item If the subsequence $(c_2,...,c_{n-1})$ contains consecutive entries that are all -2 or consecutive entries that are all two, the sequence is not minimal. In the former case, it is possible to reduce the number of entries by the transformation (\ref{eqn:surgerytwo2}). In the latter case, one can apply the inverse transformation to (\ref{eqn:surgerytwo1}) in order to reduce the number of entries.
\item If the sequence  $(c_1,...,c_n)$ contains a subsequence of the form $(a,2,3,3,...,3,3,2,b)$ or $(a,-2,-3,-3,...,-3,-3,-2,b)$, it is not minimal since then the number of entries can be reduced by a transformation of type (\ref{eqn:232}).
\item If none of the above is true, a further reduction of the sequence is not possible, i.e. the presentation is minimal. In this case, the transformations (\ref{eqn:surgery1})-(\ref{eqn:surgery3}) and (\ref{eqn:232}) are not applicable, and transformations of type (\ref{eqn:surgerytwo2}) or inverse to (\ref{eqn:surgerytwo1}) only convert entries with value -2 into entries with value 2 and vice versa without creating a new subsequence of consecutive twos or consecutive -2's.\\
Especially, a subsequence of the form $(a,2,4,2,b)$ with $a,b \neq 2$ can not be further reduced, although it is possible to reduce the entry in the middle from four down to two by changing the sign of the two 2's. But then, the new subsequence is $(a-1,-2,2,-2,b-1)$, and hence, there is no subsequence of consecutive twos or -2's. Further application of the inverse transformation to (\ref{eqn:surgerytwo1}) will only result in $(a-1,-3,-2,-3,b-1)$. The same argumentation holds for subsequences of the form $a,-2,-4,-2,b$ with $a,b \neq -2$.
\item If a presentation is minimal, it is unique up to a choice for the signs of all entries with absolute value 2 except perhaps the first and the last entry of the sequence. This comes from the fact that for $k=2$ resp. $k=1$, the transformations (\ref{eqn:surgerytwo1}) and (\ref{eqn:surgerytwo2}) do not change the number of entries but convert twos in -2's and vice versa. Since the transformations (\ref{eqn:surgery1})-(\ref{eqn:surgery3}), the inverse transformation to (\ref{eqn:surgerytwo1}) for $k>2$ and the transformation (\ref{eqn:surgerytwo2}) for $k>1$ reduce the number of entries by at least one, this is the only possibility of choice.\\
Note that flipping the sign of such an entry also modifies the two neighbouring entries by one. As another consequence, it follows that a minimal presentation of $M(c_1,...,c_n)$ of $A$ is unique if and only if all entries $c_2,...,c_{n-1}$ have an absolute value of at least 3.
\item As already noted above, this presentation is not unique if the sequence contains entries with value 2 or -2 (except perhaps for the first or the last entry of the sequence) and both neighbouring entries (if existing) are different from the entry itself.\\
The following theorem shows that introducing the convention to choose positive sign for all these entries with absolute value two, the presentation becomes unique. Hence, every matrix $A \in PSL(2,\Z)$ has a unique minimal presentation $A=M(c_1,...,\\c_k)$ with $c_2,...,c_k \notin \{-2,-1,0,1\}$ where (without considering the first and the last entry of the sequence) no two consecutive entries can be two, and, between each pair of twos, there must be at least one entry that is not three. In the following theorem, the statement will be formalized and proven formally.
\end{itemize}
\label{bla:min}
\end{bla}

\begin{sz}
Let $A=\begin{pmatrix} a&b \\ c&d \end{pmatrix} \in PSL(2,\Z)$.\\
Then, there exists exactly one presentation of $A$ that is minimal and does not contain a -2 as an entry except perhaps as first or last entry. In other words, if the presentation is denoted by $A=M(c_1,...,c_n)$, for $i \in \{2,...,n-1\}$, it is $c_i \notin \{-2,-1,0,1\}$.
\label{sz:verylong}
\end{sz}
\begin{bw}
To prove this statement, we start with the unique presentation of $A$ as in corollary \ref{cor:unice} and show that there is an algorithm for replacing all subsequences that contradict minimality as mentioned above and that this algorithm must terminate at some point. For the uniqueness of this minimal presentation, it remains to show that the resulting minimal presentation of $A$ does not depend on the order of the operations that are performed to reduce the number of entries.\\ \\
For $c=0$, it is $A=M(m,0)$ with $m=\frac{b}d=\frac{b}a$. This presentation is the only one of length two that fulfils the conditions above and is minimal. For $c \neq 0$, the proof will be split into two parts:\\ \\
(I) Existence: Consider the unique presentation
\begin{equation}
A=\begin{pmatrix} a&b \\ c&d \end{pmatrix} = M(q_1,...,q_k,m,0)
\end{equation}
with $\frac{a}c=[[q_1,...,q_k]]$ as in corollary \ref{cor:unice} and transform this presentation into one fulfilling the conditions above. First, consider $m$. If $m=0$, remove $(0,0)$ at the end. If $m=-1$, perform a transformation of type (\ref{eqn:surgery2}) to remove the -1. If $m=-2$, apply transformation (\ref{eqn:surgerytwo2}) with $k=1$ to switch the sign of $m$. If $m=1$, perform the inverse transformation to (\ref{eqn:surgerytwo3}) with the maximal value of $k$ such that the entry before the consecutive twos is at least 3 (or $q_1$) and is then still at least two after the transformation. Note that for $q_k>2$, the number of consecutive twos before $m$ is zero, and one has to choose $k=1$ in equation (\ref{eqn:surgerytwo3}), but the argumentation still holds.\\
Since $q_2,...,q_k \geq 2$, now the only issue that can contradict minimality of the sequence is that the sequence may contain a subsequence of the form $(a,2,3,...,3,2,b)$ or of the form $(a,2,...,2,b)$ (note that $a$ and $b$ are not the entries of the matrix $A$ here). To remove these subsequences, apply the following transformations:
\begin{equation}
(a,\underbrace{2,...,2}_{k~\mathrm{times}},b) \to (a-1,-(k+1),b-1)~\textrm{and}
\label{eqn:op1}
\end{equation}
\begin{equation}
(a,2,\underbrace{3,...,3}_{k~\mathrm{times}},2,b) \to (a-1,\underbrace{-3,...,-3}_{(k+1)~\mathrm{times}},b-1).
\label{eqn:op2}
\end{equation}
The initial sequence does not contain any entries in $\{-2,-1,0,1\}$ as inner entries. Hence, after applying transformation (\ref{eqn:op1}) as often as necessary, this is still the case if $k$ is always chosen maximal. If one first removes all subsequences of the form $(a,2,...,2,b)$ this way, there is no subsequence of consecutive twos left, and hence, if applying equation (\ref{eqn:op2}), it is guaranteed that either $a=c_1$ or $a \neq 2$ and analogously $b=c_n$ or $b \neq 2$, and hence, applying transformation (\ref{eqn:op2}) does also not create any inner entries with values in $\{-2,-1,0,1\}$.\\
Note that after replacing all subsequences of type $(a,2,3,..,3,2,b)$ via transformation (\ref{eqn:op2}), there might be some newly created subsequences of consecutive twos (although only such of length two), and then one has to go back to step one and apply again transformation (\ref{eqn:op1}). However, both of the named operations reduce the total number of entries, and at some point the process will end up in a sequence that is minimal. An example for a subsequence that makes it necessary to go back to step one again is $(a,2,3,3,2,3,2,b)$. After applying equation (\ref{eqn:op2}), one has $(a-1,-3,-3,-3,2,2,b)$, and hence, one more transformation (\ref{eqn:op1}) is needed to obtain $(a-1,-3,-3,-4,-3,b-1)$.\\ \\
(II) Uniqueness: So far, it has been shown that if operations of type (\ref{eqn:op1}) are always applied before operations of type (\ref{eqn:op2}) if possible, the obtained sequence fulfils the conditions above. What is still to be shown is that the obtained sequence does not depend on the order of the transformations. If the two subsequences of the above types are separated by an entry that is at least four or the subsequences are not directly neighboring to each other, this is obvious. Hence, the four following cases have to be considered:
\begin{itemize}
\item Case (1): A subsequence of the type $(a,2,3,...,3,2,3,...,3,2,b)$. Then, applying transformation (\ref{eqn:op2}) twice yields
\begin{equation}
\begin{aligned}
&(a,2,\underbrace{3,...,3}_{k~\mathrm{times}},2,\underbrace{3,3,...,3}_{l~\mathrm{times}},2,b)\\
\to ~&(a-1, \underbrace{-3,...,-3}_{(k+1)~\mathrm{times}},2,\underbrace{3,...,3}_{(l-1)~\mathrm{times}},2,b)\\
\to ~&(a-1, \underbrace{-3,...,-3}_{k~\mathrm{times}},-4,\underbrace{-3,...,-3}_{l~\mathrm{times}},b-1).
\end{aligned}
\end{equation}
Since the action of (\ref{eqn:op2}) is symmetric, replacing the second subsequence before the first one gives the same result as reverting the sequence, switching $k$ and $l$, applying the same procedure as above and reverting back the result.
\item Case (2): A subsequence of the type $(a,2,...,2,3,2,...,2,b)$. Then, applying transformation (\ref{eqn:op1}) twice yields
\begin{equation}
\begin{aligned}
&(a,\underbrace{2,...,2}_{k~\mathrm{times}},3,\underbrace{2,...,2}_{l~\mathrm{times}},b)\\
\to ~&(a-1, -(k+1),\underbrace{2,2,...,2}_{(l+1)~\mathrm{times}},b)\\
\to ~&(a-1, -(k+2),-(l+2),b-1).
\end{aligned}
\label{eqn:2and3}
\end{equation}
Again, due to the symmetry of (\ref{eqn:op1}), reverting the sequence, switching $k$ and $l$ and reverting back the result, replacing the second subsequence first does not give a different result.
\item Case (3): A subsequence of the type $(a,2,3,...,3,2,2,...,2,b)$. This is the most difficult case since there, if the subsequence containing the threes before is replaced first, there will be entries with value one or zero in some intermediate steps, and, the symmetry of the transformations (\ref{eqn:op1}) and (\ref{eqn:op2}) can also not be used. So, first check what happens when starting with the subsequence containing the threes:
\begin{equation}
\begin{aligned}
&(a,2,\underbrace{3,...,3}_{k~\mathrm{times}},\underbrace{2,...,2}_{l~\mathrm{times}},b)\\
\to ~&(a-1,\underbrace{-3,...,-3}_{(k+1)~\mathrm{times}},1,\underbrace{2,...,2}_{(l-2)~\mathrm{times}},b)\\
\to ~&(a-1,\underbrace{-3,...,-3}_{(k+1)~\mathrm{times}},0,-(l-1),b-1)\\
\to ~&(a-1,\underbrace{-3,...,-3}_{k~\mathrm{times}},-(l+2),b-1)
\end{aligned}
\label{eqn:equivalent}
\end{equation}
where the second step is applying (\ref{eqn:op1}) and the third step removing the zero by application of (\ref{eqn:surgery3}). If one elects to remove the one by transformation (\ref{eqn:surgery1}) instead of applying (\ref{eqn:op1}) as the second step, the part $(-3,1,2,...,2,b)$ takes the form $(-4,1,2,...,2,b)$ with one two less than before.\\
Now, one can either apply transformation (\ref{eqn:surgery1}) iteratively to remove all ones until reaching $b$. Checking the number of twos, note that this corresponds to applying operation (\ref{eqn:surgery1}) $l-1$ times altogether, and since $-3-(l-1) \cdot 1=-(l+2)$, the result is the same. If one decides to stop applying operation (\ref{eqn:surgery1}) after $k$ iterations and to apply operation (\ref{eqn:op1}) instead, one has now $l-k$ instead of $l$, but also $-(3+k)$ instead of the last $-3$, and hence, the result is still the same.\\ \\
Considering the whole subsequence and, as in the existence part of the proof, starting with the transformation of the second part by applying (\ref{eqn:op1}) first and then applying (\ref{eqn:op2}), one gets
\begin{equation}
\begin{aligned}
&(a,2,\underbrace{3,...,3}_{k~\mathrm{times}},\underbrace{2,...,2}_{l~\mathrm{times}},b)\\
\to ~&(a,2,\underbrace{3,...,3}_{(k-1)~\mathrm{times}},2,-(l+1),b-1)\\
\to ~&(a-1,\underbrace{-3,...,-3}_{k~\mathrm{times}},-(l+2),b-1)
\end{aligned}
\end{equation}
what is the same, hence also in this case the order of the operations does not matter.
\item Case (4): A subsequence of the type $(a,2,...,2,3,...,3,2,b)$. Due to the symmetry of all operations used in the previous case, reverting the sequence, applying the same argumentation and reverting back the result shows that also in that case it does not matter which operation is executed first.
\end{itemize}

Now, it remains to check that in the case $m \in \{-1,0,1,2\}$, it also does not matter whether one fixes that problem directly at the beginning of the process or at some point later. Clearly, if $q_k \geq 4$, one sees easily that it does not matter when the disallowed entry is removed. For $q_k=3$, removing $m$ may may decrease $q_k$ down to two, and hence might enable a transformation of type (\ref{eqn:op1}) or (\ref{eqn:op2}) if $m=1$. However, this transformation can then only be applied after removing the one and not before. Hence, the only relevant value for $q_k$ one has to consider is $q_k=2$. Therefore, consider the four cases $m=-2$, $m=-1$, $m=0$ and $m=1$.

\begin{itemize}
\item Case (1): $m=1$. If $q_k=2$, it is possible that $[[q_1,...,q_k]]$ ends on $(a,2,3,...,3,2)$ or $(a,2,...,2)$, and then there is a choice whether one wants to remove the one first or the named subsequence first.\\
(a) Case (1.1): $[[q_1,...,q_k]]$ ends on $(a,2,3,...,3,2)$. Then, the end of the sequence $(q_1,...,q_k,1,0)$ can be reduced the following two ways:
\begin{equation}
\begin{aligned}
&(a,2,\underbrace{3,...,3}_{k~\mathrm{times}},2,1,0) \to (a-1,\underbrace{-3,...,-3}_{(k+1)~\mathrm{times}},0,0) \to (a-1,\underbrace{-3,...,-3}_{(k+1)~\mathrm{times}})
\end{aligned}
\end{equation}
or
\begin{equation}
\begin{aligned}
&(a,2,\underbrace{3,...,3}_{k~\mathrm{times}},2,1,0) ~~&\to& ~~ (a,2,\underbrace{3,...,3}_{k~\mathrm{times}},1,-1)\\
\to ~&(a,2,\underbrace{3,...,3}_{(k-1)~\mathrm{times}},2,-2) ~~&\to& ~~ (a-1,\underbrace{-3,...,-3}_{k~\mathrm{times}},-3)\\
= ~&(a-1,\underbrace{-3,...,-3}_{(k+1)~\mathrm{times}}).
\end{aligned}
\end{equation}
We observe that the result is the same.\\
(b) Case (1.2): $[[q_1,...,q_k]]$ ends on $(a,2,...,2)$. Then, the end of the sequence $(q_1,...,q_k,1,0)$ can be reduced either by application of transformation (\ref{eqn:op1}) yielding
\begin{equation}
(a,\underbrace{2,...,2}_{k~\mathrm{times}},1,0) \to (a-1,-(k+1),0,0) \to (a-1,-(k+1))
\end{equation}
or by iteratively removing the ones that are created when the last one is removed, and hence, the two on the position before is reduced to one. In each iteration step, the last value in the sequence is reduced by one, and in the last step, $a$ is also reduced by one. Since one needs $k+1$ iterations, the result is also $(a-1,-(k+1))$. If one decides to stop the iteration process at a certain point and apply transformation (\ref{eqn:op1}) on the remaining subsequence of consecutive twos, the result is also the same for the same reasons as in the argumentation below equation (\ref{eqn:equivalent}).

\item Case (2): $m=0$. Obviously, it does not matter at what time the subsequence $(0,0)$ is removed without changing anything else.

\item Case (3): $m=-1$. If $q_k=2$, it is possible that $[[q_1,...,q_k]]$ ends on $(a,2,3,...,3,2)$ or $(a,2,...,2)$, and then there is a choice whether one wants to remove the -1 first or the named subsequence first.\\
(a) Case (3.1): $[[q_1,...,q_k]]$ ends on $(a,2,3,...,3,2)$. Then, if the end of the sequence $(q_1,...,q_k,-1,0)$ is reduced by applying transformation (\ref{eqn:op2}) first, one gets
\begin{equation}
\begin{aligned}
&(a,2,\underbrace{3,...,3}_{k~\mathrm{times}},2,-1,0) ~~&\to& ~~ (a-1,\underbrace{-3,...,-3}_{(k+1)~\mathrm{times}},-2,0)\\
\to ~&(a-1,\underbrace{-3,...,-3}_{k~\mathrm{times}},-2,2,1) ~~&\to& ~~ (a-1,\underbrace{-3,...,-3}_{(k-1)~\mathrm{times}},-2,2,3,1)\\
\to ~&(a-1,\underbrace{-3,...,-3}_{(k-2)~\mathrm{times}},-2,2,3,3,1) ~~&\to& ~~  ...\\
\to ~&(a-1,-2,2,\underbrace{3,...,3}_{k~\mathrm{times}},1) ~~&\to& ~~ (a,2,\underbrace{3,...,3}_{(k+1)~\mathrm{times}},1)
\end{aligned}
\end{equation}
what is also the result when directly removing the -1. The reason one has to go through this whole calculation above and cannot just stop at some point is that otherwise the resulting sequence contains an entry with value -2 contradicting the conditions of the statement (the presentation is already minimal after the first step but does not conform with the canonical choice that all entries with absolute value two in the inner part of the sequence have positive sign).\\
(b) Case (3.2): $[[q_1,...,q_k]]$ ends on $(a,2,...,2)$. Then, the end of the sequence $(q_1,...,q_k,-1,0)$ can be reduced the following two ways:
\begin{equation}
(a,\underbrace{2,...,2}_{k~\mathrm{times}},-1,0) \to (a-1,-(k+1),-2,0) \to (a-1,-k,2,1)
\end{equation}
or
\begin{equation}
(a,\underbrace{2,...,2}_{k~\mathrm{times}},-1,0) \to (a,\underbrace{2,...,2}_{(k-1)~\mathrm{times}},3,1) \to (a-1,-k,2,1).
\end{equation}
We observe that the result is the same.

\item Case (4): $m=-2$. If $q_k=2$, it is possible that $[[q_1,...,q_k]]$ ends on $(a,2,3,...,3,2)$ or $(a,2,...,2)$.\\
(a) Case (4.1): $[[q_1,...,q_k]]$ ends on $(a,2,3,...,3,2)$. Then, if the end of the sequence $(q_1,...,q_k,-1,0)$ is reduced by application of transformation (\ref{eqn:op2}) first, one obtains directly
\begin{equation}
(a,2,\underbrace{3,...,3}_{k~\mathrm{times}},2,-2,0) \to (a-1,\underbrace{-3,...,-3}_{(k+1)~\mathrm{times}},-3,0) = (a-1,\underbrace{-3,...,-3}_{(k+2)~\mathrm{times}},0).
\end{equation}
Converting the -2 into a two first yields
\begin{equation}
(a,2,\underbrace{3,...,3}_{k~\mathrm{times}},2,-2,0) \to (a,2,\underbrace{3,...,3,3}_{(k+1)~\mathrm{times}},2,1) \to (a-1,\underbrace{-3,...,-3}_{(k+2)~\mathrm{times}},0)
\end{equation}
what is the same result.\\
(b) Case (4.2): $[[q_1,...,q_k]]$ ends on $(a,2,...,2)$.  Then, if the end of the sequence $(q_1,...,q_k,-1,0)$ is reduced by application of transformation (\ref{eqn:op1}) first, one obtains directly
\begin{equation}
(a,\underbrace{2,...,2}_{k~\mathrm{times}},-2,0) \to (a-1,-(k+1),-3,0).
\end{equation}
Converting the -2 into a two first yields
\begin{equation}
\begin{aligned}
&(a,\underbrace{2,...,2,2}_{k~\mathrm{times}},-2,0) ~~&\to& ~~ (a,\underbrace{2,...,2}_{(k-1)~\mathrm{times}},3,2,1)\\
\to ~& (a-1,-k,2,2,1) ~~&\to& ~~ (a-1,-(k+1),-3,0),
\end{aligned}
\end{equation}
and again, the result is the same. Note that in this calculation, it does not depend on whether it is executed like above or the subsequence $(2,3,2)$ is replaced by $(-3,-3)$ first since this is the same situation as in case (4) of the previous part of the proof.
\end{itemize}
Now, the proof can be easily finished. It is already shown that there is only one reduction $(d_1,...,d_r)$ of the unique sequence $(q_1,...,q_k,m,0)$ such that the result is a minimal presentation of $A$ without a $-2$ as an inner entry.\\
Suppose $M(e_1,...,e_s)$ is another minimal presentation of $A$ fulfilling this property. Then, by theorem \ref{sz:sequence}, it is possible to expand the sequence $(e_1,...,e_n)$ to $(q_1,...,q_k,m,0)$ by replacing all negative entries using the corresponding operations. Inverting this process gives a method to get from $(q_1,...,q_k,m,0)$ to $(e_1,...,e_n)$ with the operations considered above, and since this reduction process is unique, it follows $(d_1,...,d_r)=(e_1,...,e_s)$.
\end{bw}

\begin{bsp}
(a) Consider the  sequence $(c_1,...,c_n)=(3,4,-5,-1,2,-3)$. Then, for
\begin{equation*}
A = M(3,4,-5,-1,2,-3) = \begin{pmatrix} 503&152\\182&55 \end{pmatrix},
\end{equation*}
one can compute $\frac{503}{182}=[[3,5,2,2,2,5,2,2]]$ and $\left \lceil \frac{55}{182} \right \rceil=1$. Hence, $A=M(3,5,2,2,2,5,2,\\2,1,0)$ in $PSL(2,\Z)$. To obtain a minimal presentation, first replace all subsequences of consecutive twos and obtain $A=M(3,4,-4,3,-3,0,0)$, and then one easily sees that $A=M(3,4,-4,3,-3)$ is a minimal presentation.\\
Since all entries are at least 3 in terms of absolute values, this minimal presentation is unique. The minimal presentation $M(3,4,-4,3,-3)$ is obtained from the given presentation $M(3,4,-5,-1,2,-3)$ by removing the entry -1 via an operation of the form (\ref{eqn:surgery2}).\\ \\
(b) Consider the matrix $B=\begin{pmatrix} -144&-55\\55&21 \end{pmatrix}$. It is $\frac{-144}{55}=[[-2,2,3,3,3,2]]$ and $\left \lceil \frac{21}{55} \right \rceil =1$. Hence,
\begin{equation*}
B=M(-2,2,3,3,3,2,1,0)=M(-3,-3,-3,-3,-3,0,0)=M(-3,-3,-3,-3,-3)
\end{equation*}
where the last expression is minimal and unique since all entries are larger at least 3 in terms of absolute values.\\ \\
(c) Consider the matrix
\begin{equation*}
C=M(4,-2,-3,-2,2,2,3,2)=\begin{pmatrix} 591&-374\\128&-81 \end{pmatrix}.
\end{equation*}
Taking the given sequence, the subsequence $(4,-2,-3,-2,2)$ can be replaced by $(5,3,3,\\3)$ applying transformation (\ref{eqn:232}) with $l=1$, and appending $(0,0)$ then results in the presentation $C=M(5,3,3,3,2,3,2,0,0)$, and one can deduce that $\frac{591}{128}=[[5,3,3,3,2,3,2]]$ since this presentation is the unique one as in corollary \ref{cor:unice}.\\
To obtain a minimal presentation, omit the $(0,0)$ at the end. Now, there is still a subsequence of the form $(2,3,2)$, but since it is at the very end of the sequence, there is no further possibility of reduction. Hence, $C=M(5,3,3,3,2,3,2)$ is a minimal presentation.\\
Instead of omitting $(0,0)$, it is also possible to transform the subsequence $(3,2,3,2,0)$ into $(2,-3,-3,-1)$ first, and then the resulting sequence is $(5,3,3,2,-3,-3,-1,0)$. Then, the -1 can be eliminated, and the sequence takes the form $C=M(5,3,3,2,-3,-2,1)$. This is another minimal presentation of $C$. Taking this one, the first one is obtained by applying (\ref{eqn:surgerytwo2}) twice (first with respect to the -2 what creates another -2 on which it is applied the second time).\\ \\
(d) Consider the matrix $D=\begin{pmatrix} -145&52\\382&-137 \end{pmatrix}$. The presentation as in corollary \ref{cor:unice} is $D=M(0,3,3,4,5,3,0,0)$. If $(0,0)$ at the end is omitted, one already has a minimal presentation $M(0,3,3,4,5,3)$. This minimal presentation is unique since except the first entry, all entries are at least 3 in terms of absolute values.\\ \\
(e) Consider the matrix $E=\begin{pmatrix} 119&-44\\46&-17 \end{pmatrix}$. The presentation as in corollary \ref{cor:unice} is $E=M(3,3,2,4,3,0,0)$. If $(0,0)$ at the end is omitted, one already has a minimal presentation $M(3,3,2,4,3)$. This minimal presentation is not unique since it contains a two in the middle. Hence, transformation (\ref{eqn:op1}) with $k=1$ can be applied with respect to that two, and the result is $M(3,2,-2,3,3)$, another minimal presentation of $E$.
\end{bsp}
\section{Minimal Presentation of Conjugacy Classes}
\label{sec:conj}
In this section, the statements made for the minimal presentation of a matrix in $PSL(2,\Z)$ are extended on conjugacy classes, i.e. some criteria are developped for a presentation of a matrix $A \in PSL(2,\Z)$ to be the shortest presentation of any matrix in the conjugacy class of $A$. In fact, it turns out that with some exceptions for special cases, the same rules as for single matrices apply with the difference that the exceptions for the beginning and end of the sequence $(c_1,...,c_k)$ fall away.\\
This is due to the fact that if a presentation of $A$ is given, any cyclic permutation of the sequence of this presentation defines another matrix in the same conjugacy class. Entries that cannot be removed without conjugation can now be moved in the middle of the sequence by a suitable conjugation and then removed.\\
As long as the length of a minimal presentation of a conjugacy class is at least two, this minimal presentation becomes unique modulo cyclic permutation of the sequence if the value -2 is disallowed for all coefficients.\\
In the first subsection, proper and pure matrices and sequences are defined. In the second subsection, it is shown that cyclic permutations of a sequence representing a matrix $A$ will represent another matrix in the same conjugacy class. The result of the third subsection is that non-proper matrices cannot have a presentation that is minimal for the whole conjugacy class if their minimal presentation has a length of at least three. In the fourth subsection, the converse is shown: Proper presentations are always minimal for the whole conjugacy class. The last subsection deals with the special case of presentations of length 2.\\
For more background information about this topic, see also section 2 of \cite{Katok} and section 5 of \cite{Zagier}. A similar investigation is performed in the sections 7.1 and 7.2 of \cite{Farey} for continuant matrices where all coefficients are positive. In chapter 2 of \cite{Katok}, the minimal presentation for a given conjugacy class of a matrix $A$ allowing only positive coefficients is derived for all matrices with a trace of absolute value larger than two.\\
In the whole section, eventual minus signs in front of matrices are skipped since the considered group is $PSL(2,\Z)$.

\subsection{Proper and Pure Sequences and Matrices}

\begin{df}
\hspace{1mm}
\begin{itemize}
\begin{em}
\item A sequence $(c_1,...,c_n)$ is called proper if it fulfils the following conditions:\\
1) It does not contain entries with value -1, 0 or 1.\\
2) It does not contain consecutive entries with value 2.\\
3) It does not contain consecutive entries with value -2.\\
4) It does not both start and end with 2 or both start and end with -2.\\
5) It does not contain subsequences of the form $(2,3,...,3,2)$ or $(-2,-3,...,-3,-2)$, and any cyclic permutation of the sequence does not contain such subsequences either. In other words, there is also no such subsequence if we consider the cyclic sequence that is generated when connecting the end of the sequence with its beginning.\\
6) It is not $(c_1,...,c_n)=(\pm 2,\pm 3,...,\pm 3)$ or any cyclic permutation of $(c_1,...,c_n)$ is equal to $(\pm 2,\pm 3,...,\pm 3)$.
\item A presentation $M(c_1,...,c_n)$ is called proper if the sequence $(c_1,...,c_n)$ is proper.
\item A matrix $A$ is called proper if there exists a proper presentation $A=M(c_1,...,c_n)$ of A, i.e. if the minimal presentations of $A$ are proper.
\item A sequence $(c_1,...,c_n)$ is called pure if it is proper and it does not contain entries with value -2. A presentation $M(c_1,...,c_n)$ is called pure if the sequence $(c_1,...,c_n)$ is pure. A matrix $A$ is called pure if it has a pure presentation.
\end{em}
\end{itemize}
\label{df:pure}
\end{df}

\begin{bsp}
(a) The sequence $(2,3,5,2,4,-6,-8,2,-4)$ is proper and pure.\\
(b) The sequence $(2,4,-4,-5,2,3,3,3)$ is not proper since $(4,-4,-5,2,3,3,3,2)$ is a cyclic permutation of it and contains $(2,3,3,3,2)$ as a subsequence.
\end{bsp}

\subsection{Conjugation as a Cyclic Permutation of the Sequence}

\begin{lm}
Let $A=M(a_1,...,a_n)$ be a presentation of $A$, $B=M(a_n)$ and $C=M(a_1)$. Then $M(a_n,a_1,...,a_{n-1})$ is a presentation of $BAB^{-1}$, and $M(a_2,...,a_n,a_1)$ is a presentation of $C^{-1}AC$.
\label{lm:presentation}
\end{lm}
\begin{bw}
Since $(a,0,-a,0)$ is a quiddity cycle for all integers $a$, i.e. $M(a,0,-a,0)=Id$ in $PSL(2,\Z)$, it is $B^{-1}=M(0,-a_n,0)$ and $C^{-1}=M(0,-a_1,0)$. Conjugating with $B$, in $PSL(2,\Z)$ one obtains
\begin{equation}
BAB^{-1}=M(a_n,a_1,...,a_n,0,-a_n,0)=M(a_n,a_1,...,a_{n-1},0,0)=M(a_n,a_1,...,a_{n-1})
\end{equation}
where for the second equality, equation (\ref{eqn:calc3}) and the fact that $a_n-a_n=0$ were used. Analogously, it is
\begin{equation}
C^{-1}AC=M(0,-a_1,0,a_1,...,a_n,a_1)=M(0,0,a_2,...,a_n,a_1)=M(a_2,...,a_n,a_1).
\end{equation}
\end{bw}

\begin{cor}
If $A=M(a_1,...,a_n) \in PSL(2,\Z)$, then all matrices $M(a_{k+1},...,a_n,a_1,\\...,a_k)$ with $0 \leq k < n$ are in the same conjugacy class as $A$. In other words, a cyclic permutation of the sequence of the entries of the presentation will give another matrix in the same conjugacy class.
\label{cor:round}
\end{cor}
\begin{bw}
Set $B_i=M(a_i)$ for $i \in \{k+1,...,n\}$. Then apply the previous lemma stepwise for $B=B_n,...,B_{k+1}$. Then, with $B=B_{k+1} \cdot ... \cdot B_n=M(a_{k+1},...,a_n)$, one obtains
\begin{equation}
BAB^{-1}=M(a_{k+1},...,a_n,a_1,...,a_n,0,-a_n,...,-a_{k+1},0)=M(a_{k+1},...,a_n,a_1,...,a_k).
\end{equation}
\end{bw}

\subsection{Finding Matrices with Shorter Presentations in the same Conjugacy Class}

\begin{lm}
If a matrix $A$ is not proper and its minimal presentations have a length of at least 3, there exists another matrix in the conjugacy class of $A$ that is proper and has a shorter minimal presentation, or that has a minimal presentation of length smaller than 3.
\label{lm:reduzione}
\end{lm}
\begin{bw}
Consider a minimal presentation $M(a_1,...,a_n)$ of $A$. Since the difference between the condition for minimality and the condition for propriety is basically that for minimality, there are some exceptions of the rules at the beginning or at the end of the sequence, it is possible to rotate the sequence by a conjugation as in lemma \ref{lm:presentation} to obtain a presentation that is not minimal. In some cases, it is also more convenient to flip the sign of an entry with absolute value 2 in order to find a matrix in the same conjugacy class with a shorter presentation. Concretely, there are the following possible cases for the minimal presentation of $A$:
\begin{itemize}
\item (1) $A=M(a_1,...,a_{n-1},\pm 1)$ or $A=M(\pm 1,a_2,...,a_n)$. Assume without loss of generality $A=M(a_1,...,a_{n-1},\pm 1)$. Otherwise, consider the matrix $A':=M(\pm1)^{-1}AM(\pm1)=M(a_2,...,a_n,\pm 1)$ which is in the same conjugacy class and of the named form.\\
Conjugate $A$ by $C=M(a_1)^{-1}$ to obtain by application of equation (\ref{eqn:calc1}) resp. (\ref{eqn:calc2})
\begin{equation}
CAC^{-1}=M(a_2,...,a_{n-1},\pm 1,a_1)=M(a_2,...,a_{n-1} \mp 1,a_1 \mp 1).
\end{equation}
\item (2) $A=M(a_1,...,a_{n-1},0)$ or $A=M(0,a_2,...,a_n)$. Assume without loss of generality $A=M(a_1,...,a_{n-1},0)$. Otherwise, consider the matrix $A':=M(0)^{-1}AM(0)=M(a_2,...,a_n,0)$  which is in the same conjugacy class and of the named form. Conjugate $A$ by $C=M(a_1)^{-1}$ to obtain by application of equation (\ref{eqn:calc3})
\begin{equation}
CAC^{-1}=M(a_2,...,a_{n-1},0,a_1)=M(a_2,...,a_{n-1}+a_1).
\end{equation}
\item (3) $A=M(a_1,...,a_{n-2},\pm 2,\pm 2)$, $A=M(\pm 2,\pm 2,a_3,...,a_n)$ or $A=M(\pm 2,a_2,...,a_{n-1},\\ \pm 2)$. Assume without loss of generality $A=M(a_1,...,a_{n-2},\pm 2,\pm 2)$. Otherwise, rotate the presentation by a suitable conjugation into one of this form.
Now, (\ref{eqn:powerone}) resp. (\ref{eqn:powertwo}) can be applied with $k=1$, and one obtains $A=M(a_1,...,a_{n-2}\mp 1,\mp 2,\pm 1)$. Then, applying the same procedure as in case (1) with $C=M(a_1)^{-1}$ yields
\begin{equation}
CAC^{-1}=M(a_2,...,a_{n-2}\mp 1,\mp 2,\pm 1,a_1)=M(a_2,...,a_{n-2}\mp 1,\mp 3,a_1\mp 1).
\end{equation}
\item (4) If the presentation of $A$ contains a subsequence $(2,3,...,3,2)$, assume that $A=M(a_1,...,a_m,2,3,...,3,2)$ ($m \leq n-3)$. Otherwise, transform the presentation by a suitable conjugation into one of this form. Now, iterative application of  equation (\ref{eqn:powerone}) with $k=1$ yields
\begin{equation}
\begin{aligned}
A&=M(a_1,...,a_m,2,\underbrace{3,...,3}_{l \textrm{ times}},2)\\
&=M(a_1,...,a_m-1,-2,2,\underbrace{3,...,3}_{(l-1)\textrm{ times}},2)\\
&=M(a_1,...,a_m,-1,-3,-2,2,\underbrace{3,...,3}_{(l-2)\textrm{ times}},2)\\
&=...=M(a_1,...,a_m-1,\underbrace{-3,...,-3}_{(l-2) \textrm{ times}},-2,2,3,2)\\
&=M(a_1,...,a_m-1,\underbrace{-3,...,-3}_{(l-1) \textrm{ times}},-2,2,2)\\
&=M(a_1,...,a_m-1,\underbrace{-3,...,-3}_{l \textrm{ times}},-2,1)
\end{aligned}
\label{eqn:sixtynine}
\end{equation}
where $l=n-m-2$ and the length of the sequence has not changed yet. Once again, the situation is the same as in case (1), and $A$ can be conjugated by $C=M(a_1)^{-1}$. Applying equation (\ref{eqn:calc1}) with $k=1$ afterwards gives
\begin{equation}
CAC^{-1}=M(a_2,...,a_m-1,\underbrace{-3,...,-3}_{l \textrm{ times}},-2,1,a_1)=M(a_2,...,a_m-1,\underbrace{-3,...,-3}_{(l+1) \textrm{ times}},a_1-1).
\end{equation}
The case of a subsequence of the form $(-2,-3,...,-3,-2)$ can be treated analogously.
\item (5) $A=M(2,3,...,3)$ (with $m$ threes) or the minimal presentation of $A$ is given by a cyclic permutation of this sequence. Then, as in the cases above, assume without loss of generality $A=M(3,2,3,...,3)$. An analogous calculation to the one in equation (\ref{eqn:sixtynine}) shows that it is also possible to write $A=M(2,-3,...,-3,-2,2)$ (with $(m-2)$ -3's) and then apply the same argumentation as in case 3. With $C=M(2,-3)^{-1}$, we finally conclude that
\begin{equation}
CAC^{-1}=M(\underbrace{-3,...,-3}_{(m-2) \textrm{ times}},-2,2,2,-3)=M(\underbrace{-3,...,,-3-3,-3}_{m \textrm{ times}},-4).
\end{equation}
The case $A=M(-2,-3,...,-3)$ is analogous.
\end{itemize}
In all cases, either the resulting presentation is proper and shorter than the initial one, or it has a length smaller than three, or one of the reduction procedures can be applied again.\\
One has to be a little bit careful about the length $n$ of the minimal presentation of $A$. In the first two cases, it is easy to notice that for $n \geq 3$, the procedure works. In the third case, for $n \geq 4$, this holds as well. For $n=3$, $a_{n-2}$ and $a_1$ are one and the same entry and there might be some trouble. However, the first part of the following lemma shows that a similar procedure still results in a matrix of the same conjugacy class with a shorter presentation and completes the proof.\\
In case (4), there is a similar situation for $m=1$, since then $a_m$ and $a_1$ are one and the same entry. This is treated by the second part of the following lemma. For $m \geq 2$, the procedure clearly works, and for $m=0$, there are no further complications since then it is $A=M(2,3,...,3,2)$, and the argumentation of case (3) is applicable.\\
The argumentation in case (5) holds clearly for all $m \geq 3$, and hence for all $n \geq 4$. For $n=3$ resp. $m=2$, the following adaption is needed: The number of -3's in the second expression of $A$ is zero, and hence, it is $A=M(2,-2,2)$ and $C=M(2,-2)^{-1}$ is needed instead of $C=M(2,-3)^{-1}$. Finally, one obtains $CAC^{-1}=M(2,2,-2)=M(1,-2,-3)$. Another conjugation by $C'=M(-3)$ gives $M(-3,1,-2)=M(-4,-3)$.
\end{bw}

\begin{lm}
(i) Let $A=M(a,\pm 2,\pm 2)$, $A=M(\pm 2,a,\pm 2)$ or $A=M(\pm 2,\pm 2,a)$. Then, the matrix $B=M(\mp 3,a\mp 2)$ is conjugate to $A$.\\
(ii) Let $A=M(a,2,\underbrace{3,...,3}_{l \textrm{ times}},2)$. Then, the matrix $B=M(\underbrace{-3,...,-3}_{(l+1) \textrm{ times}},a-2)$ is conjugate to $A$.
\label{lm:additional}
\end{lm}
\begin{bw}
(i) Without loss of generality, assume $A=M(a,\pm 2,\pm 2)$. Otherwise conjugate with a suitable matrix $D$ such that $DAD^{-1}=M(a,\pm 2,\pm 2)$. Now, application of equation (\ref{eqn:powerone}) resp. (\ref{eqn:powertwo}) with $k=1$ gives also $A=M(a\mp 1,\mp 2,\pm 1)$. With $C=M(a\mp1)^{-1}$, it follows
\begin{equation}
CAC^{-1}=M(\mp 2,\pm 1,a\mp 1)=M(\mp 3,a\mp 2)=B.
\end{equation}
This completes the proof of the previous lemma for $n=3$ in case (3).\\
(ii) As in equation (\ref{eqn:sixtynine}), it is $A=M(a-1,\underbrace{-3,...,-3}_{l \textrm{ times}},-2,1)$. Conjugation with $C=M(a-1)^{-1}$ yields
\begin{equation}
CAC^{-1}=M(\underbrace{-3,...,-3}_{l \textrm{ times}},-2,1,a-1)=M(\underbrace{-3,...,-3}_{(l+1) \textrm{ times}},a-2)=B.
\end{equation}
This completes the proof of the previous lemma for $m=1$ in case (4).
\end{bw}

\begin{bsp}
In all following examples, the given equations hold in $PSL(2,\Z)$. Eventual minus signs in front of matrices may be dropped.\\
(a) Consider the matrix $A=M(3,3,-4,0)=\begin{pmatrix}8&-35\\3&-13\end{pmatrix}$. It is not proper since its minimal presentation ends on a zero. Conjugation with $B=M(3)^{-1}$ gives
\begin{equation}
BAB^{-1}=M(0,-3,0,3,3,-4,0,3)=M(3,-1)=\begin{pmatrix}4&3\\1&1\end{pmatrix}.
\end{equation}
This matrix is still not proper but has a minimal presentation of length two, and hence, lemma \ref{lm:reduzione} is not applicable any more. However, in this case, lemma \ref{lm:dimtwo} can be used to conclude that $M(5)$ is a matrix with an even shorter presentation in the same conjugacy class (see also example \ref{bsp:dimtwo}).\\
(b) Consider the non-proper matrix $A=M(1,-2,5,-4,2)=\begin{pmatrix}150&-67\\103&-46\end{pmatrix}$ and set $B=M(2)$. Then
\begin{equation}
D:=BAB^{-1}=M(2,1,-2,5,-4,2,0,-2,0)=M(1,-3,5,-4)=\begin{pmatrix}88&21\\67&16\end{pmatrix}.
\end{equation}
One notes that this is a shorter presentation of a matrix in the same conjugacy class. However, since the first entry in this presentation is one, $D$ is not proper as well, lemma \ref{lm:reduzione} can be applied again, and conjugation with $C=M(-4)$ yields
\begin{equation}
CDC^{-1}=M(-4,1,-3,5,-4,0,4,0)=M(-5,-4,5)=\begin{pmatrix}-100&19\\21&-4\end{pmatrix}.
\end{equation}
This is a proper matrix given in a proper presentation in the same conjugacy class as $A$. Looking at definition \ref{df:pure}, one notes that the matrix $CDC^{-1}$ is even pure.\\
(c) Consider the non-proper matrix $A=M(-2,5,5,2,-2)=\begin{pmatrix}-243&-95\\110&43\end{pmatrix}$ and set $B=M(-2,5)^{-1}$. Then
\begin{equation}
BAB^{-1}=M(0,-5,2,0,-2,5,5,2,-2,-2,5)=M(5,3,3,6)=\begin{pmatrix}208&-37\\45&-8\end{pmatrix}
\end{equation}
is a proper (and pure) matrix in the same conjugacy class.\\
(d) Consider the non-proper matrix $A=M(3,2,-4,3,-2,5,2,3)=\begin{pmatrix}4132&-1687\\1619&-661\end{pmatrix}$ and set $B=M(3,2,-4)^{-1}$. Then
\begin{equation}
\begin{aligned}
BAB^{-1}&=M(0,4,-2,-3,0,3,2,-4,3,-2,5,2,3,3,2,-4)\\
&=M(3,-2,4,-3,-3,-3,-5)=\begin{pmatrix}3266&707\\947&205\end{pmatrix}.
\end{aligned}
\end{equation}
This presentation is proper but not pure. Hence, $BAB^{-1}$ is a proper matrix in the conjugacy class of $A$ with a minimal presentation of length 7. The obtained presentation can be easily transformed into the proper presentation $M(4,2,5,-3,-3,-3,5)$ by application of operation (\ref{eqn:powertwo}) with $k=1$.\\
(e) Consider the matrix $A=M(3,3,2,3,3)=\begin{pmatrix}80&-31\\31&-12\end{pmatrix}$. Conjugating with $B_1=M(3)^{-1}$ gives
\begin{equation}
B_1AB_1^{-1}=M(0,-3,0,3,3,2,3,3,3)=M(3,2,3,3,3)=M(2,-3,-3,-2,2).
\end{equation}
Now, set $B_2=M(2,-3)^{-1}$ and obtain with $B=B_2B_1$
\begin{equation}
BAB^{-1}=M(0,3,-2,0,2,-3,-3,-2,2,2,-3)=M(-3,-3,-3,-4)=\begin{pmatrix}-76&-21\\29&8\end{pmatrix}
\end{equation}
what is a proper (and pure) presentation of a matrix in the conjugacy class of $A$.
\label{bsp:dimhigh}
\end{bsp}

\subsection{Minimality of Proper Presentations for their Whole Conjugacy Class}

\begin{lm}
Let $A \in PSL(2,\Z)$ be a proper matrix, $M(a_1,...,a_n)$ a minimal (and hence proper) presentation of $A$ of length $n \geq 3$. Furthermore, let $B=M(c)$ for an arbitrary integer $c$.\\
Then, the matrix $D:=BAB^{-1}$ has no shorter minimal presentation than $A$. Furthermore, if the length of the minimal presentations is the same (i.e. the length of the minimal presentations of $D$ is also $n$), the minimal presentations of $D$ are also proper, and, if the minimal presentations of $D$ are longer than $n$, reducing the presentations of $D$ by further conjugation always leads to another proper presentation of length $n$.
\label{lm:proof}
\end{lm}
\begin{bw}
Since $(c,0,-c,0)$ is a quiddity cycle, one observes that
\begin{equation}
D=M(c)M(a_1,...,a_n)M(c)^{-1}=M(c,a_1,...,a_n,0,-c,0)=M(c,a_1,...,a_{n-1},a_n-c,0).
\end{equation}
Now, there are the following cases:
\begin{itemize}
\item (1) $c=a_n-1$:
\begin{equation}
D=M(c,a_1,...,a_n-c,0)=M(a_n-1,a_1,...,a_{n-1},1,0)=M(a_n-1,a_1,...,a_{n-1}-1,-1).
\end{equation}
If $a_{n-1} \neq 2$, this presentation of $D$ is minimal and has length $n+1$. To get rid of the -1 at the end, a conjugation as in case (1) of lemma \ref{lm:reduzione} with $C=M(a_n-1)^{-1}$ is needed, and
\begin{equation}
CDC^{-1}=M(a_1,...,a_{n-1}-1,-1,a_n-1)=M(a_1,...,a_n)=A.
\end{equation}
If $a_{n-1}=2$, it is
\begin{equation}
\begin{aligned}
D&=M(a_n-1,a_1,...,a_{n-1}-1,-1)\\
&=M(a_n-1,a_1,...,a_{n-2},1,-1)\\
&=M(a_n-1,a_1,...,a_{n-2}-1,-2).
\end{aligned}
\end{equation}
Recall that $(a_1,...,a_n)$ is proper and $a_{n-1}=2$. Hence, it follows $a_{n-2} \neq 2$, $a_n \neq 2$ and either $a_n \neq 3$ or $a_1 \neq 2$ since otherwise, there is a contradiction to point 4) or 5) of definition \ref{df:pure}. Therefore, this presentation of $D$ is minimal and proper and has length $n$.
\item (2) $c=a_n+1$:
\begin{equation}
D=M(c,a_1,...,a_n-c,0)=M(a_n+1,a_1,...,a_{n-1},-1,0)=M(a_n+1,a_1,...,a_{n-1}+1,1).
\end{equation}
If $a_{n-1} \neq -2$, this presentation of $D$ is minimal and has length $n+1$. To get rid of the 1 at the end, a conjugation as in case (1) of lemma \ref{lm:reduzione} with $C=M(a_n+1)^{-1}$ is needed, and
\begin{equation}
CDC^{-1}=M(a_1,...,a_{n-1}+1,1,a_n+1)=M(a_1,...,a_n)=A.
\end{equation}
If $a_{n-1}=-2$, it is
\begin{equation}
\begin{aligned}
D&=M(a_n+1,a_1,...,a_{n-1}+1,1)\\
&=M(a_n+1,a_1,...,a_{n-2},-1,1)\\
&=M(a_n+1,a_1,...,a_{n-2}+1,2).
\end{aligned}
\end{equation}
Recall that $(a_1,...,a_n)$ is proper and $a_{n-1}=-2$. Hence, it follows $a_{n-2} \neq -2$, $a_n \neq -2$ and either $a_n \neq -3$ or $a_1 \neq -2$ since otherwise, there is a contradiction to point 4) or 5) of definition \ref{df:pure}. Therefore, this presentation of $D$ is minimal and proper and has length $n$.
\item (3) $c=a_n-2$:
\begin{equation}
D=M(c,a_1,...,a_n-c,0)=M(a_n-2,a_1,...,a_{n-1},2,0).
\end{equation}
If $a_{n-1} \neq 2$, this presentation of $D$ is minimal and has length $n+2$. To get rid of the 0 at the end, a conjugation as in case (2) of lemma \ref{lm:reduzione} with $C=M(a_n-2)^{-1}$ is needed, and
\begin{equation}
CDC^{-1}=M(a_1,...,a_{n-1},2,0,a_n-2)=M(a_1,...,a_n)=A.
\end{equation}
If $a_{n-1}=2$, it is
\begin{equation}
D=M(a_n-2,a_1,...,a_{n-2},2,2,0)=M(a_n-2,a_1,...,a_{n-2}-1,-3,-1)
\end{equation}
Recall that $(a_1,...,a_n)$ is proper and $a_{n-1}=2$. Hence, it follows $a_{n-2} \neq 2$ and either $a_{n-2} \neq 3$ or $a_{n-3} \neq 2$ since otherwise, there is a contradiction to point 4) or 5) of definition \ref{df:pure}. Therefore, this presentation of $D$ is minimal and has length $n+1$. To get rid of the -1 at the end, a conjugation as in case (1) of lemma \ref{lm:reduzione} with $C=M(a_n-2)^{-1}$ is needed, and
\begin{equation}
\begin{aligned}
CDC^{-1}&=M(a_1,...,a_{n-2}-1,-3,-1,a_n-2)\\
&=M(a_1,...,a_{n-2}-1,-2,a_n-1)\\
&=M(a_1,...,a_{n-2},2,a_n)=M(a_1,...,a_n)=A.
\end{aligned}
\end{equation}
\item (4) $c=a_n+2$:
\begin{equation}
D=M(c,a_1,...,a_n-c,0)=M(a_n+2,a_1,...,a_{n-1},-2,0).
\end{equation}
If $a_{n-1} \neq -2$, this presentation of $D$ is minimal and has length $n+2$. To get rid of the 0 at the end, a conjugation as in case (2) of lemma \ref{lm:reduzione} with $C=M(a_n+2)^{-1}$ is needed, and
\begin{equation}
CDC^{-1}=M(a_1,...,a_{n-1},-2,0,a_n+2)=M(a_1,...,a_n)=A.
\end{equation}
If $a_{n-1}=-2$, it is
\begin{equation}
D=M(a_n+2,a_1,...,a_{n-2},-2,-2,0)=M(a_n+2,a_1,...,a_{n-2}+1,3,1)
\end{equation}
Recall that $(a_1,...,a_n)$ is proper and $a_{n-1}=-2$. Hence, it follows $a_{n-2} \neq -2$ and either $a_{n-2} \neq -3$ or $a_{n-3} \neq -2$ since otherwise, there is a contradiction to point 4) or 5) of definition \ref{df:pure}. Therefore, this presentation of $D$ is minimal and has length $n+1$. To get rid of the 1 at the end, a conjugation as in case (1) of lemma \ref{lm:reduzione} with $C=M(a_n+2)^{-1}$ is needed, and
\begin{equation}
\begin{aligned}
CDC^{-1}&=M(a_1,...,a_{n-2}+1,3,1,a_n+2)\\
&=M(a_1,...,a_{n-2}+1,2,a_n+1)\\
&=M(a_1,...,a_{n-2},-2,a_n)=M(a_1,...,a_n)=A.
\end{aligned}
\end{equation}
\item (5) $|c-a_n|>2$: Then $M(c,a_1,...,a_n-c,0)$ is a minimal presentation of $D$ and has length $n+2$. To get rid of the 0 at the end, a conjugation as in case (2) of lemma \ref{lm:reduzione} with $C=M(c)^{-1}$ is needed, and
\begin{equation}
CDC^{-1}=M(a_1,...,a_n-c,0,c)=M(a_1,...,a_n)=A.
\end{equation}
\item (6) $c=a_n$: This is precisely the case of lemma \ref{lm:presentation}, and it is
\begin{equation}
D=M(c,a_1,...,a_n-c,0)=M(a_n,a_1,...,0,0)=M(a_n,a_1,...,a_{n-1}).
\end{equation}
This presentation of $D$ is of length $n$, and its sequence is a cyclic permutation of $(a_1,...,a_n)$. Since the latter was proper, the resulting presentation of $D$ is proper as well.
\end{itemize}
\end{bw}

\begin{cor}
If a matrix $A$ is proper and its minimal presentations have at least length 3, its minimal presentations are also minimal for the conjugacy class of $A$, i.e., it is not possible to transform $A$ by conjugation into a matrix with a shorter minimal presentation.
\label{cor:proof}
\end{cor}
\begin{bw}
Suppose $A=M(a_1,...,a_n)$ is a proper presentation of $A$ and let $B \in PSL(2,\Z)$ arbitrary. Let furthermore $M(b_n,...,b_1)$ be a presentation of $B$. Then, it is
\begin{equation}
BAB^{-1}=M(b_n) \cdot ... \cdot M(b_1) A M(b_1)^{-1} \cdot ... \cdot M(b_n)^{-1}.
\end{equation}
Define furthermore for $i \in \{1,...,n\}$
\begin{equation}
A_i:=M(b_i) \cdot ... \cdot M(b_1) A M(b_1)^{-1} \cdot ... \cdot M(b_i)^{-1}.
\end{equation}
Using the previous lemma, one can now deduce iteratively that all matrices $A_i$ have no shorter minimal presentation than $A$ for $i \leq n$. Since $A_n=BAB^{-1}$, this is also valid for $BAB^{-1}$. Since $B \in PSL(2,\Z)$ was arbitrary, the statement of this corollary follows.
\end{bw}

\begin{sz}
Let $A \in PSL(2,\Z)$ and $M(a_1,...,a_n)$ with $n \geq 3$ be a minimal presentation of $A$. Then, the following is true:\\
(i) If $A$ is proper and $M(a_1,...,a_n)$ is a pure presentation of $A$, $M(a_1,...,a_n)$ and its cyclic permutations are all pure minimal presentations of matrices in the conjugacy class of $A$, and $M(a_1,...,a_n)$ is unique as a pure presentation of $A$.\\
(ii) If $A$ is proper and $M(a_1,...,a_n)$ is not a pure presentation of $A$, it is possible to construct a pure minimal presentation of a matrix in the conjugacy class of $A$ out of the given presentation, and the length of this pure minimal presentation is also $n$.\\
(iii) If $A$ is not proper, it is either possible to construct a pure matrix $B$ in the same conjugacy class taking the given minimal presentation of $A$ or to find a matrix $B$ in the conjugacy class of $A$ with a presentation of length at most 2.
\end{sz}
\begin{bw}
(i) If $M(a_1,...,a_n)$ is pure, all its cyclic permutations are also pure. Looking at the proof of lemma \ref{lm:proof}, it is clear that taking a pure presentation $M(a_1,...,a_n)$ and conjugating by $M(c)$ with $c \in \Z$, the result is a pure presentation if and only if $c=a_n$, and then, it is a cyclic permutation of the initial presentation. Applying the same argumentation as in corollary \ref{cor:proof} and case (6) of lemma \ref{lm:proof}, one can deduce that the pure presentations of matrices in the conjugacy class of $A$ are exactly the presentations of the form $M(\tilde a_1,...,\tilde a_n)$ where $(\tilde a_1,...,\tilde a_n)$ is a cyclic permutation of $(a_1,...,a_n)$.\\
The existence of two different pure presentations of the same matrix $A$ is impossible since it would contradict theorem \ref{sz:verylong}.\\
(ii) Since $M(a_1,...,a_n)$ is minimal and $A$ is proper, it follows that $M(a_1,...,a_n)$ is proper. If it is not pure, equation (\ref{eqn:powertwo}) with $k=1$ can be applied for all entries with value \mbox{-2} except perhaps the first and the last entry. If one of these entries is also -2, it is possible to move this entry away from its position by conjugation with $M(-2)^{-1}$ resp. $M(-2)$ and afterwards applying equation (\ref{eqn:powertwo}). Then, the resulting presentation is pure by construction.\\
(iii) Rotate the sequence $(a_1,...,a_n)$ by conjugation and apply the equations (\ref{eqn:powerone}), (\ref{eqn:powertwo}) and (\ref{eqn:calc1})-(\ref{eqn:calc3}) in a suitable way such that all contradictions to propriety are removed to obtain a proper presentation. The concrete procedures in the different cases are given in the proof of lemma \ref{lm:reduzione}. Then, if the resulting presentation has still length 3 or more, use part (ii) to transform this proper presentation into a pure one if necessary. Otherwise, the result is a presentation of the desired matrix $B$ with length at most 2.
\end{bw}

\begin{bsp}
In all following calculations, the given equations hold in $PSL(2,\Z)$. Eventual minus signs in front of matrices may be dropped.\\
Let $A=\begin{pmatrix}47&-224\\17&-81\end{pmatrix} \in PSL(2,\Z)$. It is $47/17=[[3,5,2,2,2]]$ and $\lceil -81/17 \rceil = -4$. Hence, according to corollary \ref{cor:unice}, it is $A=M(3,5,2,2,2,-4,0)$, and application of equation (\ref{eqn:powerone}) gives the minimal presentation $A=M(3,4,-4,-5,0)$. Since this presentation ends on a zero, $A$ is not proper, and hence also not pure. Applying lemma \ref{lm:reduzione}, one can conclude that the proper matrix
\begin{equation}
B=M(4,-4,-5,0,3)=M(4,-4,-2)=\begin{pmatrix}30&17\\7&4\end{pmatrix}
\end{equation}
belongs to the same conjugacy class.\\
Since this minimal presentation of $B$ ends on -2, it is not pure, and since this entry contradicting pureness is at the end, it cannot be transformed into a two without further conjugation. Hence, the matrix $B$ itself has no pure presentation. Therefore, $B$ is not pure. However, setting $C=M(4)^{-1}$, it is
\begin{equation}
P_1:=CBC^{-1}=M(0,-4,0,4,-4,-2,4)=M(-4,-2,4)=M(-3,2,5)=\begin{pmatrix}-32&7\\9&-2\end{pmatrix},
\end{equation}
and the given presentation of $P_1$ is pure. Hence, $P_1$ is a pure matrix in the conjugacy class of $A$. The other pure matrices in this conjugation class are $P_2=M(2,5,-3)=\begin{pmatrix}29&9\\16&5\end{pmatrix}$ and $P_3=M(5,-3,2)=\begin{pmatrix}37&-16\\7&-3\end{pmatrix}$.
\end{bsp}

\subsection{Minimal Presentations of Length Two}

\begin{lm}
Let $A=M(a_1,a_2)$ be a matrix given in its minimal presentation of length $n=2$. Then, the following statements are true:
\begin{itemize}
\item (i) The matrix $M(a_2,a_1)$ is in the same conjugacy class as $A$.
\item (ii) For $a \in \Z$, the matrices $M(a+2,2)$, $M(2,a+2)$, $M(-2,a)$ and $M(a,-2)$ are in the same conjugacy class.
\item (iii) For $a_1=\pm1$ or $a_2=\pm1$, the matrix $A$ is in the same conjugacy class as $M(a_2\mp2)$ resp. $M(a_1\mp2)$.
\item (iv) For $|a_1|,|a_2|>2$, for $|a_1|>2,a_2=0$ or $|a_2|>2,a_1=0$ and for $a_1=a_2=0$, the presentation $M(a_1,a_2)$ is minimal for the conjugacy class of $A$ and unique as a minimal presentation up to switching the two entries. Furthermore, if $|a_1|=2$ or $|a_2|=2$, the presentation is unique up to switching the two entries and the choice resulting from part (ii).
\end{itemize}
\label{lm:dimtwo}
\end{lm}
\begin{bw}
\item (i) Set $B=M(a_2)$. Then it follows
\begin{equation}
BAB^{-1}=M(a_2,a_1,a_2,0,-a_2,0)=M(a_2,a_1,0,0)=M(a_2,a_1).
\end{equation}
\item (ii) Set $B=M(1)$. Then, with $a_1=a+2$ and $a_2=2$, it follows
\begin{equation}
BAB^{-1}=M(1,a+2,2,0,-1,0)=M(1,a+2,1,0)=M(1,a+1,-1).
\end{equation}
Conjugating with $C=M(-1)$ now yields
\begin{equation}
(CB)A(CB)^{-1}=M(-1,1,a+1)=M(-2,a)
\end{equation}
And with part (i) follows that $M(a,-2)$ and $M(2,a+2)$ belong to the same conjugacy class, too.
\item (iii) Due to the symmetry of the problem (see part (i)), one may assume $a_1=\pm1$. Set furthermore $a_2=a$ and $B=M(0)$. Then it follows
\begin{equation}
BAB^{-1}=M(0,\pm1,a,0)=M(\mp1,a\mp1,0).
\end{equation}
Now, conjugating with $C=M(\mp1)^{-1}$, one obtains
\begin{equation}
(CB)A(CB)^{-1}=M(a\mp1,0,\mp1)=M(a\mp2).
\end{equation}
\item (iv) Set $B=M(c)$ for $c \in \Z$ and consider
\begin{equation}
BAB^{-1}=M(c,a_1,a_2,0,-c,0)=M(c,a_1,a_2-c,0).
\end{equation}
For $c=0$, $c=a_2$ or $c \notin \{\pm1,a_2\pm1\}$, it is clear that any combination of rotation by conjugation and operations of the types (\ref{eqn:calc1})-(\ref{eqn:calc3}) directly gives $M(a_1,a_2)$ or $M(a_2,a_1)$. For $c=\pm1$ or $c=a_2\pm1$, it may be possible to create a presentation of the kind $M(a_1-1,a_2-1,-1)$ or $M(a_1+1,a_2+1,1)$. But then, if $a_1,a_2 \neq \pm2$, the next step also gives $M(a_1,a_2)$ or $M(a_2,a_1)$. For $a_1=\pm2$ or $a_2=\pm2$, this gives precisely the proof of part (ii).
\end{bw}

\begin{cor}
A presentation $M(a_1,a_2)$ of a matrix $A \in PSL(2,\Z)$ of length two is minimal if and only if it contains no entry with an absolute value of one. It is unique up to switching the two entries if it does not contain an entry with an absolute value of two. It becomes unique up to switching the two entries if entries with value $-2$ are disallowed.
\end{cor}

\begin{bsp}
In all following calculations, the given equations hold in $PSL(2,\Z)$. Eventual minus signs in front of matrices may be dropped.\\
(a) Recall the matrix $M(3,-1)=\begin{pmatrix}4&3\\1&1\end{pmatrix}$ in example \ref{bsp:dimhigh} (a). In the following, denote $A=M(3,-1)$. With lemma \ref{lm:dimtwo}, one can now give the conjugation that transforms this matrix into $M(5)$. With $B=M(1,0)$, it follows
\begin{equation}
BAB^{-1}=M(1,0,3,-1,0,0,-1,0)=M(4,-1,-1,0)=M(5)=\begin{pmatrix}5&-1\\1&0\end{pmatrix}.
\end{equation}
(b) Consider the matrix $B_1=M(-4,-2)=\begin{pmatrix}-7&-4\\2&1\end{pmatrix}$. Applying part (ii) of lemma \ref{lm:dimtwo}, it follows that $B_2=M(-2,2)=\begin{pmatrix}-5&2\\2&-1\end{pmatrix}$  belongs to the same conjugacy class. Since in this case, part (ii) of the lemma is applicable once more, the matrix $B_3=M(2,4)=\begin{pmatrix}7&-2\\4&-1\end{pmatrix}$ must also belong to this conjugacy class.\\
(c) Applying part (ii) of the lemma on the matrix $C_1=M(7,2)=\begin{pmatrix}13&-7\\2&-1\end{pmatrix}$, it follows that the matrix $C_2=M(5,-2)=\begin{pmatrix}11&5\\2&1\end{pmatrix}$ belongs to the same conjugacy class.\\
(d) Part (iv) of the lemma guarantees that for $D_1=M(5,-3)=\begin{pmatrix}16&5\\3&1\end{pmatrix}$, the only other matrix in this conjugacy class with a minimal presentation of length two is $D_2=M(-3,5)=\begin{pmatrix}-16&3\\5&-1\end{pmatrix}$.
\label{bsp:dimtwo}
\end{bsp}

\bibliographystyle{unsrt}
\addcontentsline{toc}{section}{Bibliography}
\bibliography{Literatur}

\end{document}